\documentclass{amsart}

\usepackage{amsmath}
\usepackage{amssymb}
\usepackage{eucal}
\usepackage{amsthm}
\usepackage{verbatim}
\usepackage[hyphens]{url}
\usepackage{natbib}

\newtheorem{theorem}{Theorem}[section]
\newtheorem{proposition}[theorem]{Proposition}
\newtheorem{lemma}[theorem]{Lemma}
\newtheorem{corollary}[theorem]{Corollary}

{\theoremstyle{definition}

}

{\theoremstyle{remark}
\newtheorem{remark}[theorem]{Remark}
}

\DeclareMathOperator{\Prob}  {\mathsf{P}}
\DeclareMathOperator{\Expec} {\mathsf{E}}

\DeclareMathOperator{\rnk}{rank}
\DeclareMathOperator{\supp}{supp}

\DeclareMathOperator{\Mix}{Mix}
\DeclareMathOperator{\Cov}{Cov}

\newcommand{\msP}{ \mathnormal{\mathsf{P}} }
\newcommand{\mfF}{\mathfrak{F}}
\newcommand{\mfG}{\mathfrak{G}}
\newcommand{\mfL}{\mathfrak{L}}
\newcommand{\mfM}{\mathfrak{M}}
\newcommand{\mcA}{\mathcal{A}}
\newcommand{\mcB}{\mathcal{B}}

\newcommand{\mcF}{\mathcal{F}}
\newcommand{\mcG}{\mathcal{G}}

\newcommand{\mcL}{\mathcal{L}}

\newcommand{\mcP}{\mathcal{P}}
\newcommand{\mcX}{\mathcal{X}}

\newcommand{\mbR}{\mathbb{R}}

\newcommand{\fot}{ {\textstyle\frac{1}{2}} }

\newcommand{\impl}{~\Rightarrow~}
\newcommand{\eqval}{~\Leftrightarrow~}
\newcommand{\bydef}
    {~\overset{ \textnormal{\textsf{\tiny def}} }{=}~}

\newcommand{\stfrac}[2]{ \leavevmode\ensuremath{
    \kern.1em\raise.5ex\hbox{\footnotesize #1}
    \kern-.1em/
    \kern-.15em\lower.25ex\hbox{\footnotesize #2}} }
\newcommand{\smfrac}[2]{ \leavevmode\ensuremath{
    \kern.1em\raise.5ex\hbox{\footnotesize $#1$}
    \kern-.1em/
    \kern-.15em\lower.25ex\hbox{\footnotesize $#2$}} }
\newcommand{\Stfrac}[2]{ \leavevmode\ensuremath{
    \kern.1em\raise.5ex\hbox{#1}
    \kern-.1em/
    \kern-.15em\lower.25ex\hbox{#2}} }
\newcommand{\Smfrac}[2]{ \leavevmode\ensuremath{
    \kern.1em\raise.5ex\hbox{$#1$}
    \kern-.25em/
    \kern-.15em\lower.5ex\hbox{$#2$}} }

\begin{document} 

\title{Convergence of Estimators in LLS Analysis}

\author{Mikhail Kovtun}
\address{Mikhail Kovtun: Duke University, CDS\\
         2117 Campus Dr., Durham, NC, 27708}
\email{mkovtun@cds.duke.edu}
\thanks{The work of the first author was supported by
        the Office of Vice Provost for Research, Duke University}

\author{Anatoliy Yashin}
\address{Anatoliy Yashin: Duke University, CDS\\
         2117 Campus Dr., Durham, NC, 27708}
\email{yashin@cds.duke.edu}
\thanks{The work of the second author was supported by
        NIH/NIA grants PO1 AG08761-01, U01 AG023712 and U01 AG07198}

\author{Igor Akushevich}
\address{Igor Akushevich: Duke University, CDS\\
         2117 Campus Dr., Durham, NC, 27708}
\email{aku@cds.duke.edu}
\thanks{The work of the third author was supported by
        NIH/NIA grant P01 AG17937-05}


\begin{abstract} 
We establish necessary and sufficient conditions for consistency
of estimators of mixing distribution in linear latent structure
analysis.
\end{abstract} 

\maketitle 

\section{Introduction} 

Linear latent structure (LLS) analysis is aimed to derive
properties of population as whole and properties of individuals
from a large number of categorical measurements
made on each individual in a sample.
An exposition of LLS analysis is given in \cite{Kovtun:2005c}.

LLS analysis searches for representation of the observed
joint distribution of random variables (representing measurements)
as a mixture of {\em independent distributions},
i.e. distributions, in which random variables are mutually
independent. Such approach is common for all branches
of latent structure analysis.
The specific LLS assumption is that the mixing distribution
is supported by a low-dimensional linear subspace of the
space of independent distributions.

When dimensionality of the supporting subspace is sufficiently
smaller than the number of random variables
and under some regularity conditions,
a set of low-order moments (of order up to number of measurements)
of mixing distribution is identifiable.
Increasing the size of the sample does not increase
the number of identifiable moments---only precision of their
estimates increases.
Thus, the mixing distribution in LLS analysis
is only {\em partially identifiable}.
In \cite{Kovtun:2005c} we suggested an estimator
of the mixing distribution. This estimator, however,
does not converge to the true mixing distribution
when the sample size tends to infinity.
The natural question arise:
what is the value of such partially identifiable model
and in what sense the suggested estimator is useful?

In this article, we give one possible answer to this question.
Assume that one has an infinite sequence of possible measurements;
a particular survey uses finite number of measurements
from this sequence.
Then a sequence of estimates of mixing distribution
obtained in a sequence of surveys with increasing number
of measurements converges to the true mixing distribution.
The necessary and sufficient condition for this convergence
is pairwise orthogonality of independent distributions
being mixed.

The crucial part of the proof is the equivalence between
pairwise orthogonality of independent distribution
and orthogonality in aggregate (see lemma \ref{lm:SingSets}).
This fact is of interest on its own;
we devote a special discussion to its consequences
with respect to structure of tail $\sigma$-algebras and
the strong law of large numbers.

The fact that pairwise orthogonality of a family of measures
being mixed is necessary and sufficient for equality 
$G = \Expec(G | \mfF^\infty)$ $\nu$-a.s. (see definitions below)
was proved by A. Yashin 
for the case of countable families of measures
\citep[see][where this fact was used to prove the strong consistency
of statistical estimates in the adaptive control schemes]
{Yashin:1982, Yashin:1986}.
A generalization of this fact to the case of families
of measures of cardinality of continuum, proved by M. Kovtun,
provides a technical basis for the results of this paper.

The rest of the article is organized as follows.
We start with extensive preliminaries in order to make
the article self-sufficient for wider audience.
Section \ref{sec:TheProblem} provides all necessary
definitions and establishes basic facts.
Section \ref{sec:MainTheorem} contains a proof
of the main theorem. 
Proofs of technical facts are given in section
\ref{sec:AuxilaryLemmas}.
The article is concluded by discussion of connections 
between our proofs and structure of tail $\sigma$-algebras and
the strong law of large numbers.


\section{Preliminaries} 
\label{sec:Preliminaries}

\subsection{Basic notions and notation} 
\label{subsec:BasicNotionsAndNotation}

For a topological space $\Omega$, $\mcB(\Omega)$ denotes
its Borel $\sigma$-algebra (i.e., $\sigma$-algebra, generated
by open subsets of $\Omega$).

Let $\mu$ be a Borel probabilistic measure on a topological 
space $\Omega$.
We say that  $\mu$  {\em is carried by} a set $A \in \mcB(\Omega)$,
if $\mu(A)=1$.
A measure $\mu$ {\em is supported by} a set $A$, if $\mu$ is carried
by $A$ and $A$ is closed.
A closed set $S$ is {\em a support} of measure $\mu$,
if $S$ supports $\mu$ and $S = \cap \{A \mid A \text{ supports } \mu\}$.
The support of measure $\mu$ is denoted by $\supp(\mu)$.
Every measure on topological space with countable base (in particular,
on separable metrizable space) has a support.

Note that a set that carries a measure $\mu$ can be ``significantly
smaller'' than its support. For example, let $\{r_n\}_{n=1}^\infty$ 
be any enumeration of rational numbers of the interval $[0,1]$,
and let $\mu_r$ be a counting measure defined by $\mu_r(r_n)=2^{-n}$.
Then the set of rationals carries $\mu_r$, but the support of $\mu_r$
is the whole interval $[0,1]$.

Two measures $\mu$ and $\mu'$ are said to be {\em orthogonal}
(or {\em singular}), denoted $\mu \perp \mu'$, if there exists
a set $A \in \mcB(\Omega)$ such that $A$ carries $\mu$ and
$\Omega \setminus A$ carries $\mu'$.
For example, $\mu_r$ is orthogonal to Lebesgue measure $\lambda$
on $[0,1]$, as $\mu_r$ is carried by rationals and $\lambda$ is
carried by irrationals.

We write $\mu \ll \mu'$ if $\mu$ is absolutely continuous with
respect to $\mu'$, and $\mu \sim \mu'$ if $\mu$ is equivalent to
$\mu'$ (i.e., $\mu \ll \mu'$ and $\mu' \ll \mu$).

For a measurable space $(\Omega,\mcF)$, a set $B \in \mcF$ is called
an {\em atom of $\sigma$-algebra $\mcF$}, if for every $B' \in \mcF$,
$B' \subseteq B$, either $B'=B$ or $B'=\varnothing$.
For a space with measure $(\Omega,\mcF,\mu)$, a set $B \in \mcF$
is called an {\em atom of measure $\mu$},
if $\mu(B) > 0$ and for every $B' \in \mcF$, $B' \subseteq B$,
either $\mu(B')=\mu(B)$ or $\mu(B')=0$.
For example, atoms of $\mcB([0,1])$ are one-points sets,
while the Lebesgue measure $\lambda$ on $([0,1],\mcB([0,1]))$
has no atoms.
We say that a measure $\mu$ is {\em continuous} if it has no atoms.

We use $\mcP(\Omega)$ to denote a space of probabilistic Borel measures
on a topological space $\Omega$. We always consider $\mcP(\Omega)$
with topology of weak convergence.

In a topological space, $[B]$ denotes the closure of a set $B$.
In a metric space, $U_\varepsilon x$ denotes 
$\varepsilon$-neighborhood of point $x$,
and $U_\varepsilon B$ denotes $\varepsilon$-neighborhood of set $B$.
$I_B(\cdot)$ denotes a characteristic function of a set $B$.
For a set $B$ with $\mu(B)>0$, $\mu(\cdot|B)$ denotes a probabilistic
measure conditional on $B$ 
(i.e., $\mu(\cdot|B)(B')=\mu(B'|B)=\mu(B' \cap B)/\mu(B)$).
$B_n \downarrow B$ abbreviates the statement
``$B_1 \supseteq B_2 \supseteq \dots$ and $\bigcap_n B_n = B$.''
$B_n \downdownarrows B$ abbreviates the statement
``$B_n \downarrow B$ and 
$\forall \, \varepsilon>0 \,\,
\exists \, m \,\, \forall \, n>m : B_n \subseteq U_\varepsilon B$.''
$\delta_x$ denotes a measure concentrated at a single point $x$.

We shall use the following properties of Borel measures on 
a metric spaces:
\begin{enumerate}
\item
    $x \in \supp(\mu) \eqval \forall \varepsilon>0 :
    \mu(U_\varepsilon x)>0$
\item
    For every $x \in \supp(\mu)$ there exists a sequence of sets
    $\{B_n\}_n$ such that $\mu(B_n)>0$ and $B_n \downdownarrows \{x\}$
\item
    $\mu' \ll \mu \impl \supp(\mu') \subseteq \supp(\mu)$
\item
    If $\mu(B)>0$, then $\mu(\cdot|B) \ll \mu$ and
    $\supp(\mu(\cdot|B)) \subseteq \supp(\mu)$
\item
    If $B_n \downdownarrows \{x\}$ and $\mu_n$ is carried by $B_n$,
    then $\mu_n \xrightarrow{w} \delta_x$
\end{enumerate}

\subsection{Lebesgue spaces} 
\label{subsec:LebesgueSpaces}

The notion of Lebesgue space was introduced and its main properties
were established in \cite{Rokhlin:1949}
(English translation---\cite{Rokhlin:1952}).

Let $(\Omega,\mcF,\mu)$ be a space with measure.
In this subsection, all spaces with measure are assumed to be
complete, i.e. $B \in \mcF$, $\mu(B)=0$, and $B' \subseteq B$
imply $B' \in \mcF$.
We also assume that $\mu(\Omega)=1$.

A countable system of measurable sets $\mfG=\{\Gamma_i\}_{i=1}^\infty$
is called a {\em basis} of $(\Omega,\mcF,\mu)$, if it possesses
the following two properties:
\begin{enumerate}
\item[($\mfL$)]
    $\forall\, B \in \mcF \,\,\,\exists\, C \in \sigma(\mfG) :
    B \subseteq C \wedge \mu(B)=\mu(C)$
\item[($\mfM$)]
    $\forall\, \omega \neq \omega' \in \Omega \,\,\,\exists\, i :
    (\omega  \in \Gamma_i \wedge \omega' \not\in \Gamma_i) \vee
    (\omega' \in \Gamma_i \wedge \omega  \not\in \Gamma_i)$
\end{enumerate}
A space with measure is called {\em separable}, if it has a basis.

If $\Omega$ is a second-countable $T_0$ topological space
(in particular, if $\Omega$ is a separable metric space)
and $\mu$ is a measure defined on its Borel $\sigma$-algebra
$\mcB(\Omega)$, then $(\Omega,\mcB_\mu(\Omega),\bar{\mu})$,
where $\mcB_\mu(\Omega)$ is Lebesgue completion of Borel 
$\sigma$-algebra with respect to measure $\mu$ and $\bar{\mu}$
is continuation of $\mu$ onto $\mcB_\mu(\Omega)$,
is a separable space with measure,
and topology base serves as a basis in the above sense.

Let $(\Omega,\mcF,\mu)$ be a separable space with measure,
and let $\mfG=\{\Gamma_i\}_i$ be its basis.
Consider intersections
\begin{equation}   
\label{eq:RokhlinIntersections}
\bigcap_{i=1}^\infty \Upsilon_i, \qquad
\text{where either $\Upsilon_i=\Gamma_i$ or 
                   $\Upsilon_i=\Omega\setminus\Gamma_i$}
\end{equation}   

Due to property ($\mfM$), every such intersection cannot contain
more than one point. Furthermore, every one-point set is representable
in form (\ref{eq:RokhlinIntersections})
(to obtain a set consisting of $\omega$, take
$\Upsilon_i=\Gamma_i$ if $\omega \in \Gamma_i$, and
$\Upsilon_i=\Omega\setminus\Gamma_i$ otherwise).

A separable space with measure $(\Omega,\mcF,\mu)$ is said to be
{\em complete with respect to basis} $\mfG=\{\Gamma_i\}_i$,
if every intersection (\ref{eq:RokhlinIntersections}) is nonempty
(and thus, consists of one point).
A separable space with measure $(\Omega,\mcF,\mu)$ is said to be
{\em complete$\pmod{0}$ with respect to basis} $\mfG=\{\Gamma_i\}_i$,
if it is isomorphic$\pmod{0}$ to a space with measure
$(\Omega',\mcF',\mu'$),
which is complete with respect to its basis $\mfG'=\{\Gamma'_i\}_i$,
and this isomorphism transforms base $\mfG$ to basis $\mfG'$.

Important fact is \citep[][\S2, n$^\circ$2]{Rokhlin:1949,Rokhlin:1952}: 
if a separable space with measure is complete$\pmod{0}$
with respect to some basis, then it is complete$\pmod{0}$
with respect to every other basis.

A separable space with measure is called {\em Lebesgue space},
if it is complete$\pmod{0}$ with respect to its bases.

It happens that many important spaces with measures are
Lebesgue spaces. In particular,
if $\Omega$ is a separable (in topological sense) complete
metric space and $\mu$ is a measure defined on its Borel
$\sigma$-algebra $\mcB(\Omega)$,
then $(\Omega,\mcB_\mu(\Omega),\bar{\mu})$
is Lebesgue space \citep[][\S2, n$^\circ$7]{Rokhlin:1949,Rokhlin:1952}.


\subsection{Metric structures} 
\label{subsec:MetricStructures}

Let $(\Omega,\mcF,\mu)$ be a space with measure.
Then $\mcF$ is a Boolean algebra, and thus a ring
(with respect to additive operation $\bigtriangleup$,
symmetric difference, and multiplicative operation $\cap$).
A family of $\mu$-negligible sets
\begin{equation*}   
I_\mu = \{ B \in \mcF \mid \mu(B)=0 \}
\end{equation*}   
is a $\sigma$-ideal of $\mcF$.

The quotient algebra $\smfrac{\mcF}{I_\mu}$ is called
a {\em metric structure} of the space $(\Omega,\mcF,\mu)$.
The Boolean algebra $\smfrac{\mcF}{I_\mu}$ is always complete
\citep[][Chapter 2]{Vladimirov:2002}.
If $\mcF_\mu$ is Lebesgue completion of $\sigma$-algebra $\mcF$
and $\bar{\mu}$ is the corresponding continuation of measure $\mu$,
then $\smfrac{ \mcF_\mu }{ I_{\bar{\mu}} }$ is isomorphic
(as a Boolean algebra) to $\smfrac{\mcF}{I_\mu}$.

If $B$ and $B'$ belong to the same coset with respect to ideal $I_\mu$,
then $\mu(B)=\mu(B')$; thus, $\mu$ can be considered as being
defined on $\smfrac{\mcF}{I_\mu}$.
The quotient algebra $\smfrac{\mcF}{I_\mu}$ can be made a metric
space by defining a distance between two cosets $a$ and $b$
as $\rho(a,b) \bydef \mu(a \bigtriangleup b)$.
$(\smfrac{\mcF}{I_\mu},\rho)$ is always a complete metric space,
and it is separable in topological sense if and only if
$(\Omega,\mcF,\mu)$ is separable in Rokhlin's sense.

Two Lebesgue spaces are isomorphic$\pmod{0}$ as spaces with measure
if and only if their metric structures are isomorphic as
metric spaces \citep[][\S2, n$^\circ$6]{Rokhlin:1949,Rokhlin:1952}.

For more details on metric structures, see
\cite{Vladimirov:2002} and \cite{Rokhlin:1949,Rokhlin:1952}.


\subsection{Sufficient $\sigma$-algebras} 
\label{subsec:Sufficient}

In this article, we use a notion of sufficient $\sigma$-algebra,
which is slightly more general than a notion of sufficient statistics
and better servers for our purposes.

Let $\{ \msP_\theta \}_{\theta \in \Theta}$ be a family of
probabilistic measures defined on measurable space $(\Omega,\mcF)$.
A $\sigma$-algebra $\mcG \subseteq \mcF$ is said to be
{\em sufficient for family} $\{ \msP_\theta \}_{\theta \in \Theta}$,
if there exists a function $P(\omega,B)$, defined for all
$\omega \in \Omega$ and $B \in \mcF$, such that:
\begin{enumerate}
\item
    For every fixed $B \in \mcF$, $P(\omega,B)$ is $\mcG$-measurable
    (as function of $\omega$).
\item
    For every fixed $\omega \in \Omega$, $P(\omega,B)$ is
    a probabilistic measure on $(\Omega,\mcF)$.
\item
    For every fixed $B \in \mcF$,
    $P(\omega,B) = \msP_\theta(B | \mcG)(\omega)$
    $\msP_\theta$-a.s.
    (i.e., $P(\omega,B)$ as function of $\omega$ is a variant
    of conditional probability $\msP_\theta(B|\mcG)$).
\end{enumerate}

In this terminology, a statistic $T(\omega)$ is sufficient
for family $\{ \msP_\theta \}_{\theta \in \Theta}$,
if $\sigma$-algebra $\sigma(T)$, generated by random variable
$T(\omega)$, is sufficient.



\section{The problem} 
\label{sec:TheProblem}

We consider an infinite sequence of random variables
$\{X_j\}_j$; variable $X_j$ takes values in a finite set
$\{1,\dots,L_j\}$.
We consider these variables as being defined on a probabilistic space
$\mcA = \prod_{j=1}^\infty \{1,\dots,L_j\}$;
$X_j(a) = X_j(a_1,\dots) = a_j$.
The space $\mcA$, endowed with Tikhonov topology, is compact
and metrizable by metric $\rho(a,a') = 1/\inf\{j \mid a_j \neq a'_j\}$;
it is complete with respect to this metric.
The Borel $\sigma$-algebra $\mcB(\mcA)$ on $\mcA$ coincides 
with the $\sigma$-algebra,
generated by random variables $\{X_j\}_j$;
thus, joint distributions of $\{X_j\}_j$ are described by Borel
measures on $\mcA$. We use $\mcP(\mcA)$ to denote the space of
all probabilistic measures on $\mcA$.

We always consider $\mcP(\mcA)$ with topology of weak convergence.
Topology of weak convergence on $\mcP(\mcA)$ is metrizable,
$\mcP(\mcA)$ is compact, and thus separable
(\citealp[][IV.6.3]{Dunford:1958};
\citealp[][IV.3.4.4,5]{Kolmogorov:1972}).

Among all joint distributions of $\{X_j\}_j$, 
we distinguish {\em independent} ones,
i.e. those distributions, in which random variables $\{X_j\}_j$
are mutually independent.
To specify an independent distribution, one needs to specify only
probabilities $\beta_{jl} = \Prob(X_j=l)$.
Then, due to independence, one has for cylinders
$\Prob(X_{j_1}=l_1 \wedge \dots \wedge X_{j_p}=l_p) =
\beta_{j_1 l_1} \cdots \beta_{j_p l_p}$,
and as cylinders compose a topology base, this uniquely extends
to $\mcB(\mcA)$. 
Thus, any independent distribution is uniquely described by an
infinite-dimensional vector 
$\beta = (\beta_{11},\dots,\beta_{1 L_1},\dots,
          \beta_{j1},\dots,\beta_{j L_j},\dots)$.
To specify a distribution, such vector must satisfy conditions:
\begin{equation}   
\label{eq:betaCond}
\begin{cases}
0 \le \beta_{jl} \le 1          & \text{for all $j$ and $l$} \\
\sum_{l=1}^{L_j} \beta_{jl} = 1 & \text{for all $j$}
\end{cases}
\end{equation}   

The set of vectors satisfying (\ref{eq:betaCond}) is a convex
infinite-dimensional body $P$ in $\mbR^\infty$.
Let $\msP_\beta$ denote the independent distribution over $A$
corresponding to a vector $\beta \in P$.

Our goal is to investigate a class of mixtures of independent
distributions, i.e. those distributions $\msP$ on $\mcA$ which
can be represented in form:
\begin{equation}   
\label{eq:genMix}
\msP(A) = \int_P \msP_\beta(A) \,\mu(d\beta) 
\qquad \text{for all } A \in \mcB(\mcA)
\end{equation}   

\noindent
where $\mu$ is a probabilistic measure on $P$.
For this definition of mixture to be correct, one needs to show
that the mapping $\beta \mapsto \msP_\beta(A)$ is measurable
for every Borel $A \subseteq \mcA$.
We show a stronger fact:

\begin{proposition}   
\label{pr:Cont1}
The mapping $\beta \mapsto \msP_\beta$ is continuous with respect
to Tikhonov topology on $\mbR^\infty$ and topology of weak convergence
on $\mcP(\mcA)$.
\end{proposition}

\begin{proof}
It is sufficient to show that for every converging sequence
$\beta^n \rightarrow \beta$ in $P$, the sequence $\msP_{\beta^n}$
weakly converges and $\lim \msP_{\beta^n} = \msP_\beta$.
For this, in turn, one needs to show that 
$\msP_{\beta^n}(C) \rightarrow \msP_\beta(C)$ for every cylinder
$C$ with $\msP_\beta(\partial C)=0$
\citep[][III.1.5]{Shiryaev:2004}.
In $\mcA$ every cylinder is open-closed set; thus, this must
be shown for all cylinders.
But for a cylinder  $C$ with base $\{a_{j_1},\dots,a_{j_p}\}$ one has
$\msP_{\beta^n}(C)=\prod_q \beta^n_{j_q a_{j_q}}$, 
and required convergence follows directly from the convergence 
$\beta^n \rightarrow \beta$.
\end{proof}   

\begin{corollary}   
\label{cr:Meas1}
For every $A \in \mcB(\mcA)$, the mapping $\beta \mapsto \msP_\beta(A)$
is measurable
(w.r.t. corresponding Borel $\sigma$-algebras).
\end{corollary}

\begin{proof}
By \citet[][III.55,60]{Dellacherie:1978}, the mapping
$\msP \mapsto \msP(B)$ is measurable.
Thus, our mapping is measurable as a composition of two measurable
mappings.
\end{proof}   

We abbreviate the equation (\ref{eq:genMix}) to 
$\msP = \Mix(\mu)$.
We also use $\msP_\mu$ to denote $\Mix(\mu)$.

One of the questions that we are interested in is conditions
for identifiability of mixtures of independent distributions
\citep[in sense of][]{Teicher:1960},
i.e. under what conditions the mixture $\Mix(\mu)$ uniquely
defines the mixing measure $\mu$.
It is easy to show, however, that without additional restrictions
{\em any} distribution on $\mcA$ can be represented as a mixture
of independent ones, and (except degenerate cases) every distribution
has infinitely many such representations.

For a distribution $\mu$ on $P$, let $\mcL(\mu)$ denote the smallest
linear subspace of $\mbR^\infty$ supporting $\mu$ (this subspace
is a linear span of $\supp(\mu)$).
We say that a mixing distribution $\mu$ has rank $K$, $\rnk(\mu)=K$,
if $\dim(\mcL(\mu))=K$.
We say that a distribution $\msP \in \mcP(\mcA)$ has rank
not greater than $K$, $\rnk(\msP) \le K$, if $\msP = \Mix(\mu)$
for some measure $\mu$ with $\rnk(\mu) = K$;
further, $\rnk(\msP)=K$, if $\rnk(\msP) \le K$ 
and $\rnk(\msP) \nleq K-1$.

The question of identifiability can be now splitted into two
subquestions:

\begin{itemize}
\item
    {\em For a distribution of rank $K$, what are conditions 
    for identifiability of a $K$-dimensional linear subspace
    carrying a mixing measure?
    More precisely: what are conditions for implication
    \begin{equation*}
    \left. \begin{array}{l}
      \rnk(\msP)=K \\
      \msP=\Mix(\mu)=\Mix(\mu') \\
      \rnk(\mu)=\rnk(\mu')=K
    \end{array} \right\}  \impl  \mcL(\mu)=\mcL(\mu')
    \end{equation*}
    to be true?}
\item
    {\em For a fixed $K$-dimensional subspace $\mcL$ of $\mbR^\infty$,
    what are conditions for identifiability of mixing measure 
    carried by this subspace?
    More precisely: what are conditions for implication
    \begin{equation*}
    \left. \begin{array}{l}
      \Mix(\mu)=\Mix(\mu') \\
      \mcL(\mu)=\mcL(\mu')=\mcL
    \end{array} \right\}  \impl  \mu=\mu'
    \end{equation*}
    to be true?}
\end{itemize}

We shall address these questions full detail in another article.
Here, we mention a sufficient condition of identifiability
established (in slightly different form) in \cite{Kovtun:2005c}.
Let $C = \left( \Cov_\mu(\beta_{jl},\beta_{j'l'}) \right)_{jl,j'l'}$
be (infinite)
covariance matrix of the mixing measure $\mu$.
Suppose that all minors of $C$ of size $K+1$ are degenerate
and all minors of size $K$ are non-degenerate.
Then $\rnk(\mu)=K$ and any mixing measure of rank not greater than $K$,
which produces the same mixture, coincides with $\mu$.

When the mixing measure is identifiable, the question arise
how the mixing measure can be estimated.
The main topic of the present article is the investigation of
properties of the estimator, outlined in 
\cite{Kovtun:2005c}.
The formal description of this estimator is given below.

Let us fix a $K$-dimensional subspace $\mcL$ of $\mbR^\infty$
such that $P_0 = P \cap \mcL$ is a $(K-1)$-dimensional body.
Let us also fix $\lambda^1, \dots, \lambda^K$, a set of $K$ 
linearly independent vectors from $P_0$.
Then every vector $\beta \in P_0$ can be uniquely represented
in form $\beta = \sum_k g_k \lambda^k$, and $\sum_k g_k = 1$.
Thus, we have a linear isomorphism $\beta : Q \rightarrow P_0$,
where $Q$ is a $(K-1)$-dimensional body in $\mbR^K$, belonging 
to the hyperplane $\sum_k g_k = 1$.
For brevity, we write $\msP_g$ instead of $\msP_{\beta(g)}$. 

Because the mapping $\beta : Q \rightarrow P_0$ is obviously continuous,
we immediately obtain from proposition \ref{pr:Cont1}:

\begin{proposition}   
\label{pr:Cont2}
The mapping $g \mapsto \msP_g$ is continuous with respect
to standard Euclidean topology on $Q$ and topology of weak convergence
on $\mcP(\mcA)$.
\end{proposition}   

\begin{proposition}   
\label{pr:Meas2}
For every $B \in \mcB(\mcA)$,
the mapping $g \mapsto \msP_g(B)$ is measurable
(w.r.t. corresponding Borel $\sigma$-algebras).
\end{proposition}   

Due to the isomorphism $\beta : Q \rightarrow P_0$, every measure
$\mu_{P_0}$ on $P_0$ induces measure $\mu_Q$ on $Q$
(by letting $\mu_Q(B) = \mu_{P_0}(\beta(B))$) and vice versa.
We usually drop the index and write just $\mu$, as it is
obvious from the context what measure is considered.

Thus, every probabilistic measure $\mu$ on $Q$ generates
a probabilistic measure $\msP_\mu$ by means of (\ref{eq:genMix}).

Now we develop a slightly different view on the problem.
Measure $\mu$ may be considered as a distribution of random variable
$G$ that takes values in $Q$.
Random variables $G, X_1, \dots$ may be considered as defined
on a common probabilistic space $\Omega = Q \times \mcA$;
$G(g,a_1,\dots)=g$ and $X_j(g,a_1,\dots)=a_j$.

A distribution $\mu$ of $G$ induces a joint distribution $\nu$
of $G, X_1, \dots$: for Borel sets $B \in \mcB(Q)$ and 
$A \in \mcB(\mcA)$,
\begin{equation}   
\label{eq:nuDef}
\nu(B \times A) = \int_B \msP_g(A) \,\mu(dg)
\end{equation}   

When $A$ is a cylinder with base $\{a_{j_1},\dots,a_{j_p}\}$,
one obtains
\begin{equation*}
\nu(B \times A) = 
\int_B \prod_{q=1}^p \sum_{k=1}^K g_k \lambda^k_{j_q a_{j_q}} \,\mu(dg)
\end{equation*}
By the Kolmogorov's theorem, the measure $\nu$ uniquely extends
from sets of form $B \times C$ to the Borel $\sigma$-algebra 
of $\Omega$.

The marginal distribution of $\nu$ on $\mcA$ is $\msP_\mu$,
and the marginal distribution of $\nu$ on $Q$ is $\mu$.
Note that $\nu$ {\em is not} a product of $\mu$ and $\msP_\mu$.

Let $\mcF^n = \sigma(X_1,\dots,X_n)$, an increasing sequence of 
$\sigma$-subalgebras of $\mcB(\Omega)$,
and let $\mcF^\infty = \sigma(X_1,\dots) = \sigma(\bigcup_n \mcF^n)$.
Consider conditional expectations $\Expec(G \mid \mcF^n)$
and $\Expec(G \mid \mcF^\infty)$.
By L\'evi theorem, 
$\Expec(G \mid \mcF^n) \rightarrow \Expec(G \mid \mcF^\infty)$,
$\nu$-a.s. and in the sense of $L^1$.

On the other hand, every element of $\mcF^n$ and $\mcF^\infty$
has a form $Q \times A$, $A \in \mcB(\mcA)$
(in other words, sets $Q \times \{a\}$ are {\em atomic} w.r.t.
$\sigma$-algebras $\mcF^n$, $\mcF^\infty$).
Thus, values $\Expec(G \mid \mcF^n)(g,a)$, 
$\Expec(G \mid \mcF^\infty)(g,a)$ do not depend on $g$.
This allows us to define functions $e_n : \mcA \rightarrow Q$,
$e_\infty : \mcA \rightarrow Q$ by letting 
$e_n(a)=\Expec(G \mid \mcF^n)(g,a)$ and similarly for $e_\infty$.
Again, $e_n \rightarrow e_\infty$ $\msP_\mu$-a.s. 
and in the sense of $L^1$.

Functions $e_n$, $e_\infty$ transform measure $\msP_\mu$ on $\mcA$
to measures $\hat{\mu}_n$, $\hat{\mu}_\infty$ on $Q$
by letting $\hat{\mu}_n(B) = \msP_\mu(e_n^{-1}(B))$ for every
$B \in \mcB(Q)$, and similarly for $\hat{\mu}_\infty$.
Simple, but important fact is:

\begin{proposition}   
\label{pr:muhatConv}
The sequence $\{\hat{\mu}_n\}_n$ weakly converges to $\hat{\mu}_\infty$.
\end{proposition}

\begin{proof}
It is sufficient to show \citep[][III.1.3.V]{Shiryaev:2004}
that $\int \phi(g) \,\hat{\mu}_n(dg)$
converges to $\int \phi(g) \,\hat{\mu}_\infty(dg)$ for every
uniformly continuous function $\phi : Q \rightarrow \mbR$.
By changing variables, one obtains
\begin{align*}
\int \phi(g) \,\hat{\mu}_n(dg)      &=
                    \int \phi(e_n(a)) \,\msP_\mu(da) \\
\int \phi(g) \,\hat{\mu}_\infty(dg) &=
                    \int \phi(e_\infty(a)) \,\msP_\mu(da)
\end{align*}
and required convergence is equivalent to
$\phi \circ e_n \xrightarrow{L^1} \phi \circ e_\infty$,
which follows from the convergence $e_n \xrightarrow{L^1}e_\infty$
due to the uniform continuity of $\phi$.
\end{proof}   

Functions $e_n$ do not depend on ``tails'' of their
arguments: if $a_j=a'_j$ for $j=1,\dots,n$, then $e_n(a)=e_n(a')$.
This allows to construct $\hat{\mu}_n$ only from knowledge
of distribution of $X_1,\dots,X_n$ and from $e_n$.
One possible way to estimate $e_n$ was established in
\cite{Kovtun:2005, Kovtun:2005c};
\cite{Kovtun:2005c} suggested to use $\hat{\mu}_n$
as an estimator of $\mu$.
Conditions for consistency of this estimator are
the main topic of the present article.

Note that $\nu$, $e_n$, $\hat{\mu}_n$, etc. depend on 
the mixing measure $\mu$.
To stress such dependency, we use sometime notation $\nu[\mu]$, 
$e_n[\mu]$, $\hat{\mu}_n[\mu]$, etc.
We shall drop $[\mu]$ whenever the mixing measure $\mu$
is obvious from the context.

We shall also use $\sigma$-algebras 
$\mcF^\infty_n = \sigma(X_n,X_{n+1},\dots)$ and the tail
$\sigma$-algebra $\mcX = \bigcap_n \mcF^\infty_n$.
These $\sigma$-algebras are considered as algebras of subsets
of either $\mcA$ or $\Omega = Q \times \mcA$ (in the latest case,
they consist of cylinders built on subsets of $\mcA$).
It should be always clear from the context what case is considered.

\section{Main theorem} 
\label{sec:MainTheorem}

The main theorem, proved in this section,
establishes necessary and sufficient conditions for
convergence $\hat{\mu}_n \xrightarrow{w} \mu$.
It is formulated as nine equivalent conditions,
which are naturally splitted into three groups.
The first group consists of conditions
(\ref{mti:converges}) and (\ref{mti:equals});
they are just statements of convergence.
The second group consists of conditions
(\ref{mti:graph}), (\ref{mti:condexp}), (\ref{mti:rcd}),
and (\ref{mti:sufficient});
these conditions are expressed in terms of random
variable $G$ and related entities.
The last group consists of conditions
(\ref{mti:singular}), (\ref{mti:Hellinger1}),
and (\ref{mti:Hellinger2}).
Condition (\ref{mti:singular}) is just statement
of pairwise orthogonality of independent measures being mixed,
while other two are criteria of orthogonality
in terms of Hellinger integrals.

The condition (\ref{mti:condexp}),
$G=\Expec(G \mid \mcF^\infty)$ $\nu$-a.s.,
looks surprising for the first glance:
$G$ as function of $\omega=(g,a)$ does not depend on $a$,
while $\Expec(G \mid \mcF^\infty)$ does not depend on $g$.
The only points where these two functions coincide
is the graph of function $e_\infty$.
But, according to condition (\ref{mti:graph}), this graph
has $\nu$-measure 1, and thus equality holds $\nu$-a.s.

This condition has an important consequence:
if the sequence of outcomes $a=(a_1,\dots)$ is generated
by an independent distribution $\mcP_g$, then $e_\infty(a)$
with probability 1 gives this distribution.
This fact justifies usage of conditional expectation
$\Expec(G \mid X_1=a_1,\dots,X_j=a_j)$ as estimation
of individual position in phase space $Q$,
suggested in \cite{Kovtun:2005c}.

\begin{theorem}
The following conditions are equivalent:
\begin{enumerate}
\item \label{mti:converges}
    For every $\mu \in \mcP(Q)$, $\hat{\mu}_n[\mu] \xrightarrow{w} \mu$
\item \label{mti:equals}
    For every $\mu \in \mcP(Q)$, $\hat{\mu}_\infty[\mu] = \mu$
\item \label{mti:graph}
    For every $\mu \in \mcP(Q)$, $\nu[\mu](\Gamma e_\infty[\mu]) = 1$
    (where $\Gamma f$ denotes the graph of a function $f$)
\item \label{mti:condexp}
    For every $\mu \in \mcP(Q)$, 
    $G = \Expec(G \mid \mcF^\infty)$ $\nu[\mu]$-a.s.
\item \label{mti:rcd}
    For every $\mu \in \mcP(Q)$, there exists a variant of regular 
    conditional distribution of $G$ with respect to $\mcF^\infty$
    (i.e., a function $P(\omega,B)$, $\omega \in \Omega$, 
    $B \in \mcB(Q)$)
    such that for every $\omega \in \Omega$, the distribution
    $P(\omega,\cdot)$ is concentrated at a single point.
\item \label{mti:sufficient}
    For every $\mu \in \mcP(Q)$, there exists a set $Q_1 \subseteq Q$
    with $\mu(Q_1)=1$ such that
    tail $\sigma$-algebra $\mcX$ is sufficient for
    the family of measures $\{ \msP_g \mid g \in Q_1 \}$.
\item \label{mti:singular}
    $\forall g', g'' \in Q :
        g' \neq g'' \impl \msP_{g'} \perp \msP_{g''}$
\item \label{mti:Hellinger1}
    $\forall g', g'' \in Q : g' \neq g'' \impl
        \prod_{j=1}^\infty \sum_{l=1}^{L_j}
        \sqrt{\beta_{jl}(g') \beta_{jl}(g'')} = 0$
\item \label{mti:Hellinger2}
    $\forall g', g'' \in Q : g' \neq g'' \impl
        \sum_{j=1}^\infty \sum_{l=1}^{L_j}
        \left( \sqrt{\beta_{jl}(g')} - \sqrt{\beta_{jl}(g'')} \right)^2
        = \infty$
\end{enumerate}
\end{theorem}

\begin{proof}
In order to make the proof of main theorem more clear,
we moved proofs of all technical facts to the next section.

\medskip

(\ref{mti:converges})$\eqval$(\ref{mti:equals}).
Follows from the proposition \ref{pr:muhatConv}.

\medskip

(\ref{mti:equals})$\impl$(\ref{mti:singular}).
Take arbitrary $g',g'' \in Q$ and take $\mu \in \mcP(Q)$ being 
concentrated at points $g',g''$ with weights $\alpha$ and $1-\alpha$,
where $0 < \alpha < 1$.
Then for every $B \in \mcB(Q)$,
$\mu(B) = \hat{\mu}_\infty(B) = \msP_\mu(e^{-1}_\infty(B))$.
In particular, for $A' = e^{-1}_\infty(\{g'\})$, one obtains
$\msP_\mu(A') = \mu(\{g'\}) = \alpha$.

Further, one obtains:
\begin{multline*}   
g' \cdot \alpha =
g' \cdot \int_{A'} \msP_\mu(da) =
\int_{A'} g' \,\msP_\mu(da) =
\int_{A'} e_\infty(a) \,\msP_\mu(da) =
\\
\int_{Q \times A'} \Expec(G \mid \mcF^\infty)(\omega) \,\nu(d\omega) =
\int_{Q \times A'} G(\omega) \,\nu(d\omega) \overset{\text{\tiny(*)}}{=}
\\
g' \cdot \nu(\{g'\} \times A') + g'' \cdot \nu(\{g''\} \times A')
    \overset{\text{\tiny(**)}}{=}
g' \cdot \alpha \msP_{g'}(A') + g'' \cdot (1-\alpha) \msP_{g''}(A')
\end{multline*}   
(equality (*) holds as $G(\omega)$ is a simple function,
taking with non-zero probability only values $g'$ and $g''$;
equality (**) follows from definition of $\nu$,
equation (\ref{eq:nuDef})). 

As $g'$ and $g''$ are linearly independent vectors, the above
equation implies that $\msP_{g'}(A') = 1$ and $\msP_{g''}(A') = 0$,
which means $\msP_{g'} \perp \msP_{g''}$, q.e.d.

\medskip

(\ref{mti:singular})$\eqval$(\ref{mti:sufficient}).
Follows from lemmas \ref{lm:TailSufficient}
and \ref{lm:InvTailSufficient}.

\medskip

(\ref{mti:singular})$\impl$(\ref{mti:rcd}).
We shall show a stronger fact, that there exists a variant of regular
conditional distribution with respect to the tail
$\sigma$-algebra $\mcX$, which has the required property
and which serves as regular conditional distribution
w.r.t. $\mcF^\infty$ as well.
Take $P((g,a),B) = \delta_{\varphi(a)}(B)$,
where $\varphi$ defined by lemma \ref{lm:InvSingSets}.
For every $\omega = (g,a)$, $P(\omega,B)$ is, by definition,
a probabilistic measure on $(Q,\mcB(Q))$, and this measure
is concentrated at a single point.
Thus, one has to show
that for every $B \in \mcB(Q)$,
$\nu(G^{-1}(B) | \mcX)(\omega) =
\Expec(I_{G^{-1}(B)} | \mcX)(\omega) = P(\omega,B)$
$\nu$-a.s.
As by lemma \ref{lm:InvSingSets} $P(\omega,B)$ is $\mcX$-measurable,
one needs to show that for every $Q \times A \in \mcX$
\begin{equation}   
\label{eq:rcdCond}
\int_{Q \times A} I_{G^{-1}(B)}(\omega) \,\nu(d\omega) =
\int_{Q \times A} P(\omega,B) \,\nu(d\omega)
\end{equation}   
Using Robbins' theorem (section \ref{subsec:PairwiseOrthogonality}),
the left-hand side is reduced to:
\begin{multline*}   
\int_{Q \times A} I_{G^{-1}(B)}(\omega) \,\nu(d\omega) =
\int_{Q \times A} I_B(g) \,\nu(d\omega) =
\\
\int_Q \left( \int_A I_B(g) \,\msP_g(da) \right) \,\mu(dg) =
\int_Q I_B(g) \left( \int_\mcA I_A(a) \,\msP_g(da) \right) \,\mu(dg) =
\\
\int_Q I_B(g) \msP_g(A) \,\mu(dg) 
\end{multline*}   
where $h$ is a homomorphism defined in section
\ref{subsec:PairwiseOrthogonality}.
The right-hand side, in turn, is reduced to:
\begin{multline*}   
\int_{Q \times A} P(\omega,B) \,\nu(d\omega) =
\int_{Q \times A} \delta_{\varphi(a)}(B) \,\nu(d\omega) =
\\
\int_Q\left(\int_A\delta_{\varphi(a)}(B)\,\msP_g(da)\right)\,\mu(dg) =
\int_Q \delta_g(B) \msP_g(A) \,\mu(dg) =
\\
\int_Q I_B(g) \msP_g(A) \,\mu(dg)
\end{multline*}   
which proves the equality (\ref{eq:rcdCond}).
Finally, note that the proof of the equality (\ref{eq:rcdCond})
remain valid, if $A$ is taken from $\mcF^\infty$.
Thus, $P(\omega,B)$ is also a regular conditional distribution
w.r.t. $\mcF^\infty$.

\medskip

(\ref{mti:rcd})$\impl$(\ref{mti:condexp}).
Let $P(\omega,B)$ be a regular conditional distribution of $G$
w.r.t. $\mcF^\infty$, such that for every $\omega$ the distribution
$P(\omega,\cdot)$ is concentrated at a single point.
Note that $P(\omega,B)$ as function of $\omega=(g,a)$ does not
depend on $g$ (as $\mcF^\infty$ consists of cylinders $Q \times A$);
this allows us to write $P(a,B)$.
Let further $\bar{g}(a)$ be the point, which the distribution
$P(a,\cdot)$ is concentrated at.
By \citet[][II.7.7]{Shiryaev:2004}, conditional expectation
can be calculated using regular conditional distribution:
\begin{equation*}   
\Expec(G | \mcF^\infty)(g',a') =
\int_Q g P(a',dg) =
\bar{g}(a')
\end{equation*}   

For arbitrary $C \in \mcB(\Omega)$ one has
$\nu(C) = \int_\mcA P(a,C_a) \msP_\mu(da)$, where
$C_a = \{ g \mid (g,a) \in C \}$.
Applying this formula to
$C_0 = \{ \omega \mid G(\omega) = \Expec(G|\mcX)(\omega) \} =
\{(g,a) \mid g = \bar{g}(a) \}$
one obtains
\begin{equation*}   
\nu(C_0) =
\int_\mcA P(a,\{\bar{g}(a)\}) \,\msP_\mu(da) =
\int_\mcA 1 \cdot \,\msP_\mu(da) = 1
\end{equation*}   
which proves that $G(\omega) = \Expec(G|\mcX)(\omega)$ $\nu$-a.s. 

\medskip

(\ref{mti:condexp})$\impl$(\ref{mti:graph}).
Follows from:
\begin{equation*}   
\Gamma e_\infty = \{ (g,a) \mid g = e_\infty(a) \} =
\{ (g,a) \mid G(g,a) = \Expec(G | \mcX)(g,a) \}
\end{equation*}   

\medskip

(\ref{mti:graph})$\impl$(\ref{mti:equals}).
First, note that for every $B \in \mcB(Q)$ one has
$(B \times \mcA) \cap \Gamma e_\infty = 
(Q \times e_\infty^{-1}(B)) \cap \Gamma e_\infty$.
Then (as $\msP(A)=1 \impl \msP(B)=\msP(B \cap A)$ and
$\mu$ and $\msP_\mu$ are marginals of $\nu$) one obtains
$\mu(B) = \nu(B\times\mcA) = \nu((B\times\mcA) \cap \Gamma e_\infty)$
and $\hat{\mu}_\infty(B) = \msP_\mu(e_\infty^{-1}(B)) =
\nu(Q \times e_\infty^{-1}(B)) =
\nu((Q \times e_\infty^{-1}(B)) \cap \Gamma e_\infty)$,
which gives $\mu(B) = \hat{\mu}_\infty(B)$ for every Borel $B$, q.e.d.

\medskip

(\ref{mti:singular})$\eqval$(\ref{mti:Hellinger1}).
By \citet[][III.9.3, theorem 3]{Shiryaev:2004}, 
$\msP_{g'} \perp \msP_{g''}$ is equivalent to
$H(\msP_{g'},\msP_{g''})=0$,
where $H(\msP,\msP')$ is a Hellinger integral of order $\frac{1}{2}$.
When random variables $X_1,\dots$ are independent
(which is the case for distributions $\msP_{g'}$ and 
$\msP_{g''}$), their joint distribution is the product of
individual distributions, and Hellinger integral for the joint
distributions is the product of Hellinger integrals for the
individual distributions \citep[][III.9.2]{Shiryaev:2004}.
The direct computation (see section \ref{subsec:HellingerIntegralsAnd})
gives
$H(\msP_{g'},\msP_{g''}) =
\prod_{j=1}^\infty 
\sum_{l=1}^{L_j} \sqrt{\beta_{jl}(g')\beta_{jl}(g'')}$,
q.e.d.

\medskip

(\ref{mti:Hellinger1})$\eqval$(\ref{mti:Hellinger2}).
Follows from lemma \ref{lm:ProdEqSum}.
\end{proof}


\section{Auxilary lemmas} 
\label{sec:AuxilaryLemmas}

\subsection{Hellinger integrals and related series} 
\label{subsec:HellingerIntegralsAnd}

Let $\msP^j_g$ denote the distribution of $X_j$ that corresponds to $g$.
$\msP^j_g$ is defined on a finite probabilistic space $\{1,\dots,L_j\}$,
and is fully determined by vector
$(\beta_{j1}(g),\dots,\beta_{j L_j}(g))$.
The Hellinger integral of order $\fot$ 
\citep[see][III.9.3]{Shiryaev:2004} 
for $\msP^j_{g'}$ and $\msP^j_{g''}$ is:
\begin{equation}   
\label{eq:HellingerInt1}
H(\msP^j_{g'},\msP^j_{g''}) = 
\sum_{l=1}^{L_j}  \sqrt{\beta_{jl}(g')\beta_{jl}(g'')}
\end{equation}   

As random variables $\{X_j\}_j$ are mutually independent in
distribution $\msP_g$, the measure $\msP_g$ is a product of measures
$\msP^j_g$;
thus, the Hellinger integral for distributions $\msP_{g'}$
and $\msP_{g''}$ is
\begin{equation}   
\label{eq:HellingerInt2}
H(\msP_{g'},\msP_{g''}) =
\prod_{j=1}^\infty H(\msP^j_{g'},\msP^j_{g''}) =
\prod_{j=1}^\infty 
    \sum_{l=1}^{L_j} \sqrt{\beta_{jl}(g')\beta_{jl}(g'')}
\end{equation}   

The important properties of the Hellinger integral are:
\begin{equation}   
\label{eq:HellingerProperties}
\text{(a)} \quad \msP' \perp \msP'' \eqval H(\msP',\msP'')=0 \qquad
\text{(b)} \quad \msP' \sim  \msP'' \impl  H(\msP',\msP'')>0
\end{equation}   

We are also interested in related quantity
\begin{equation}   
\label{eq:HelSum}
H^+(\msP_{g'},\msP_{g''}) =
    \sum_{j=1}^\infty \sum_{l=1}^{L_j}
        \left( \sqrt{\beta_{jl}(g')} - \sqrt{\beta_{jl}(g'')} \right)^2
\end{equation}   

To simplify notation, we shall use
\begin{equation}   
\label{eq:AbbrHellinger}
  H(g',g'') \bydef   H(\msP_{g'},\msP_{g''}); \qquad
H^+(g',g'') \bydef H^+(\msP_{g'},\msP_{g''})
\end{equation}   

The equality of the Hellinger integral to zero may come from
two sources: first, one of the factors in the product equals to $0$,
and second, all factors are positive but converge to zero
sufficiently quickly.

In the second case, one can employ the following considerations.
For arbitrary sequence $\{x_n\}_n$ with $0 < x_n \le 1$, one has
\begin{equation}   
\label{eq:equivalence}
\prod_n x_n = 0  \eqval
\sum_n \ln(x_n) = -\infty   \eqval
\sum_n \ln\frac{1}{x_n} = \infty    \eqval
\sum_n (1 - x_n) = \infty
\end{equation}   

\begin{remark}   
\label{rm:ImplToLeft}
If in (\ref{eq:equivalence}) $x_n$ are allowed to be zero,
all equivalences except the last one remain true; for the last
step, only implication from right to left holds.
\end{remark}   

Further, for every $j$
\begin{multline*}   
\sum_l \left( \sqrt{\beta_{jl}(g')} - \sqrt{\beta_{jl}(g'')} \right)^2 =
\\
\sum_l \beta_{jl}(g') + \sum_l \beta_{jl}(g'') - 
2 \sum_l \sqrt{\beta_{jl}(g')\beta_{jl}(g'')} =
\\
2 \left( 1 - \sum_l \sqrt{\beta_{jl}(g')\beta_{jl}(g'')} \right)
\end{multline*}   

Thus, one obtains:

\begin{lemma}   
\label{lm:HellingerSum1}
If for every $j$, $\msP^j_{g'} \sim \msP^j_{g''}$, then
\begin{equation}   
\label{eq:ProdEqSum}
H(g',g'') = 0  \eqval  H^+(g',g'') = \infty
\end{equation}   
\end{lemma}   

In the first case, however, the equivalence (\ref{eq:ProdEqSum})
is not necessarily true.
If, for example, $\beta(g') = (1,0,\frac{1}{2},\frac{1}{2},\dots)$
and $\beta(g'') = (0,1,\frac{1}{2},\frac{1}{2},\dots)$
(all random variables $X_j$ are assumed binary), one has
\begin{equation*}   
H(g',g'') = 
\left( \sqrt{1\cdot0}+\sqrt{0\cdot1} \right) \cdot
\left( \sqrt{\fot\cdot\fot}+\sqrt{\fot\cdot\fot} \right)
\cdot \dots = 
0 \cdot 1 \cdot \dots = 0
\end{equation*}   
while the right-hand series in (\ref{eq:ProdEqSum}) equals to
\begin{multline*}   
\left( \left( \sqrt{1}-\sqrt{0} \right)^2 +
       \left( \sqrt{0}-\sqrt{1} \right)^2 \right) +
\left( \left( \sqrt{\fot}-\sqrt{\fot} \right)^2 +
       \left( \sqrt{\fot}-\sqrt{\fot} \right)^2 \right) + \dots \\ 
= 2 + 0 + \dots = 2
\end{multline*}   

Nevertheless, the equality to zero of $H(g',g'')$ for {\em every}
pair $g', g''$ is equivalent to divergence $H^+(g',g'')$
for {\em every} pair $g', g''$.
We shall establish this weaker (but sufficient for
our purposes) fact in the lemma \ref{lm:ProdEqSum} below;
but before this we need a couple of additional notions and lemmas.

For $g',g'' \in Q$, let $[g',g'']$ denote a closed interval
with endpoints $g'$ and $g''$ (i.e., 
$[g',g''] = \{ (1-t) \cdot g' + t \cdot g'' \mid 0 \le t \le 1 \}$).
Similarly, $(g',g'')$ denotes an open interval and
$\langle g',g'' \rangle$ denotes the whole line passing through
$g'$ and $g''$.
Note that $Q$ is convex, bounded and closed;
thus, for every $g',g'' \in Q : [g',g''] \subseteq Q$
and intersection $\langle g',g'' \rangle \cap Q$ is a closed interval.

\begin{lemma}   
\label{lm:PointsOutside}
If $H^+(g',g'')=\infty$ and $[g',g''] \subseteq [\bar{g}',\bar{g}'']$,
then $H^+(\bar{g}',\bar{g}'')=\infty$.
\end{lemma}

\begin{proof}
Due to linearity of mapping $\beta$, one has for every $j$, $l$:
$|\sqrt{\beta_{jl}(g')}-\sqrt{\beta_{jl}(g'')}| \le
|\sqrt{\beta_{jl}(\bar{g}')}-\sqrt{\beta_{jl}(\bar{g}'')}|$.
\end{proof}   

\begin{lemma}   
\label{lm:simpleIneq}
For every $a,b \in [0,1]$ one has
\begin{equation*}
\left(\frac{3}{2}-\sqrt{2}\right)
    \cdot \left( \sqrt{a} - \sqrt{b} \right)^2
\le \left( \sqrt{\frac{a+b}{2}} - \sqrt{b} \right)^2
\end{equation*}
\end{lemma}

\begin{proof}
If $b=0$, the inequality is verified directly.
If $b>0$, one can divide both sides by $b$ and consider a function
\begin{equation*}
\phi(t)=\left(\frac{3}{2}-\sqrt{2}\right)
    \cdot \left( \sqrt{t} - 1 \right)^2
- \left( \sqrt{\frac{t+1}{2}} - 1 \right)^2
\end{equation*}
By analyzing this function, one obtains that it has a single maximum
on $(0,\infty)$ at point $1$, which equals $0$.
\end{proof}   

\begin{corollary}   
\label{cr:DivMidPoint}
$H^+(g',g'')=\infty$ $\impl$
$H^+(\frac{g'+g''}{2},g'')=H^+(g',\frac{g'+g''}{2})=\infty$.
\end{corollary}   

\begin{lemma}   
\label{lm:PairImpliesLine}
If for some $\bar{g}',\bar{g}'' \in Q$ one has 
$H^+(\bar{g}',\bar{g}'')=\infty$, then
$H^+(g',g'')=\infty$ for every 
$g',g'' \in \langle \bar{g}',\bar{g}'' \rangle \cap Q$. 
\end{lemma}

\begin{proof}
If $[\bar{g}',\bar{g}''] \subseteq [g',g'']$, the statement
follows from lemma \ref{lm:PointsOutside}.
Otherwise, by taking $g'_{(1)}=\bar{g}'$, $g''_{(1)}=\bar{g}''$
and either 
$g'_{(i+1)}=g'_{(i)}$, $g''_{(i+1)}=\frac{g'_{(i)}+g''_{(i)}}{2}$
or
$g'_{(i+1)}=\frac{g'_{(i)}+g''_{(i)}}{2}$,  $g''_{(i+1)}=g''_{(i)}$,
one comes to $g'_{(n)}$, $g''_{(n)}$ such that
$[g'_{(n)},g''_{(n)}] \subseteq [g',g'']$.
Then the required statement follows from corollary
\ref{cr:DivMidPoint} and lemma \ref{lm:PointsOutside}.
\end{proof}   

\begin{lemma}   
\label{lm:InteriorImpl}
If $\bar{g}',\bar{g}'' \in Q$, then for every 
$g',g'' \in (\bar{g}',\bar{g}'')$ one has
$\msP^j_{g'} \sim \msP^j_{g''}$ for every $j$.
\end{lemma}

\begin{proof}
The only reason for $\msP^j_{g'} \not\sim \msP^j_{g''}$
is that for some $l$, $\beta_{jl}(g')=0$ and $\beta_{jl}(g'')>0$
(or vice versa). But due to linearity of mapping $\beta$
this may happen only for endpoints of the interval
$\langle \bar{g}',\bar{g}'' \rangle \cap Q$.
\end{proof}   

\begin{corollary}   
\label{cr:InteriorImpl}
Let $\bar{g}',\bar{g}'' \in Q$. 
Then $H(g',g'')=0 \eqval H^+(g',g'')=\infty$
for every $g',g'' \in (\bar{g}',\bar{g}'')$.
\end{corollary}   

\begin{lemma}   
\label{lm:ProdEqSum}
$\forall g',g'' \in Q : H(g',g'')=0
\eqval
\forall g',g'' \in Q : H^+(g',g'')=\infty$.
\end{lemma}

\begin{proof}
The implication from right to left follows from remark
\ref{rm:ImplToLeft}.
To prove implication from left to right, assume premise and
take arbitrary $g',g'' \in Q$.
Further, take arbitrary $\bar{g}',\bar{g}'' \in (g',g'')$.
By assumption, $H(\bar{g}',\bar{g}'')=0$, and by corollary
\ref{cr:InteriorImpl} $H^+(\bar{g}',\bar{g}'')=\infty$.
Then by lemma \ref{lm:PairImpliesLine} $H^+(g',g'')=\infty$, q.e.d.
\end{proof}   


\subsection{Tail space} 
\label{subsec:TailSpace}

In this and subsequent subsections we consider random variables
$X_j$ as defined on measurable space $(\mcA,\mcB(\mcA))$;
thus,the $\sigma$-algebras $\mcF^n = \sigma(X_1,\dots,X_n)$,
$\mcF^\infty = \sigma(\bigcup_n \mcF^n)$,
$\mcF^\infty_n = \sigma(X_n,X_{n+1},\dots)$, and
$\mcX = \bigcap_n \mcF^\infty_n)$
are subalgebras of $\mcB(\mcA)$.

Consider a relation $\sim$ on $\mcA$ defined as:
\begin{equation}   
\label{eq:TailEqDef}
a' \sim a'' \bydef \exists\,n \forall\,m>n : a'_m = a''_m
\end{equation}   

Obviously, this relation is an equivalence relation;
let $\widetilde{a}$ denote a coset of $a \in \mcA$,
and $\widetilde{\mcA}$ denote a factorspace of $\mcA$.

Unfortunately, the factor topology of $\widetilde{\mcA}$
is trivial (a closure of every one-point set is the whole space),
and thus does not possess any useful for our purpose property.

One can easily verify that for every $A \in \mcX$ the implication
$a \in A \impl \widetilde{a} \subseteq A$ holds,
and that for every $a \in \mcA$ the coset $\widetilde{a}$ belongs
to the tail $\sigma$-algebra $\mcX$.
This makes the tail $\sigma$-algebra $\mcX$ isomorphic to its
image  
$\widetilde{\mcX} = \{ \widetilde{A} \subseteq \widetilde{\mcA} \mid
\exists A \in \mcX : \widetilde{A}=\{\widetilde{a} \mid a \in A\} \}$
under factor mapping.
Any measure $\msP$ on $(\mcA,\mcB(\mcA))$ can be restricted to
$\sigma$-algebra $\mcX$; this restriction gives rise to
measure $\widetilde{\msP}$ on $(\widetilde{\mcA},\widetilde{\mcX})$.
Also, any $\mcX$-measurable function $f$ on $\mcA$ 
corresponds to a $\widetilde{\mcX}$-measurable function $\widetilde{f}$
on $\widetilde{\mcA}$ defined as 
$\widetilde{f}(\widetilde{a}) = f(a)$.

The measurable space $(\widetilde{\mcA},\widetilde{\mcX})$ is called 
a {\em tail space}.
The slight advantage of considering the tail space 
instead of $(\mcA,\mcX)$ is that atoms of $\widehat{\mcX}$
(in the sense of Boolean algebra) are one-point sets,
while atoms of $\mcX$ are cosets of relation $\sim$.

In the subsequent subsections, we drop the $\sim$ sign
and use $\mcA$, $\mcX$, etc. to denote $\widetilde{\mcA}$,
$\widetilde{\mcX}$, etc.
This should not led to confusion, as any statement regarding
$(\widetilde{\mcA},\widetilde{\mcX})$ can be reformulated
as a statement regarding $(\mcA,\mcX)$, and in the most cases
such reformulation consists just of dropping $\sim$ sign.


\subsection{Orthogonality 
            of measures restricted to $\mcX$}
\label{subsec:OrthogonalityWithRespect}

For every measures $\msP'$, $\msP''$ on $(\mcA,\mcB(\mcA))$
the absolute continuity $\msP' \ll \msP''$ implies the 
absolute continuity of measures restricted to $\mcX$,
$\msP'|_\mcX \ll \msP''|_\mcX$.
We shall show that a similar implication holds for orthogonality
in the case of family of measures $\{\msP_g\}_{g \in Q}$. Namely:

\begin{lemma}   
\label{lm:OrtWrtX}
Suppose that 
$\forall g',g'' \in Q : g' \neq g'' \impl \msP_{g'} \perp \msP_{g''}$. 
Then 
$\forall g',g'' \in Q : 
g' \neq g'' \impl \msP_{g'}|_\mcX \perp \msP_{g''}|_\mcX$,
i.e. there exist $A \in \mcX$
such that $\msP_{g'}(A)=1$ and $\msP_{g''}(A)=0$.
\end{lemma}

\begin{proof}
Assume premise of the lemma.
By lemma \ref{lm:ProdEqSum},
\begin{equation*}
\textstyle
\sum_{j=1}^\infty \sum_l 
\left( \sqrt{\beta_{jl}(g')} - \sqrt{\beta_{jl}(g'')} \right)^2
\end{equation*}
diverges for every $g' \neq g''$.
This implies that for every $n$ the sum
\begin{equation*}
\textstyle
\sum_{j=n}^\infty \sum_l 
\left( \sqrt{\beta_{jl}(g')} - \sqrt{\beta_{jl}(g'')} \right)^2
\end{equation*}
also diverges.
From this, by applying considerations of subsection
\ref{subsec:HellingerIntegralsAnd} to the space
$\mcA^{(n)} = \prod_{j=n}^\infty \{1,\dots,L_j\}$,
one obtains $\msP^{(n)}_{g'} \perp \msP^{(n)}_{g''}$
for arbitrary $g' \neq g''$
(here $\msP^{(n)}_g$ denotes a measure on $\mcA^{(n)}$
defined in the way similar to the definition of $\msP_g$ on $\mcA$).
But $\msP^{(n)}_{g'}$ and $\msP^{(n)}_{g''}$ are marginals
of $\msP_{g'}$ and $\msP_{g''}$, respectively;
thus, for every $n$ there exists a set $A_n \in \mcF^\infty_n$
such that $\msP_{g'}(A_n)=1$ and $\msP_{g''}(A_n)=0$.

Now take $A = \bigcap_{n=1}^\infty \bigcup_{j=n}^\infty A_n$.
It is easy to verify that $A \in \mcX$ and $\msP_{g'}(A)=1$,
$\msP_{g''}(A)=0$, which proves the lemma.
\end{proof}   

\begin{remark}  
The condition that the orthogonality takes place for {\em every}
pair $g',g'' \in Q$ is important.
To obtain a counterexample, consider the case:
$Q = [0,1]$, $K=2$, all variables $X_1,\dots$ are binary,
$\lambda^1 = (1,0,\frac{1}{2},\frac{1}{2},\dots)$,
$\lambda^2 = (0,1,\frac{1}{2},\frac{1}{2},\dots)$, and
$\beta(g) = g \lambda^1 + (1-g) \lambda^2$.
Then $\msP_0 \perp \msP_1$, but $\msP_0|_\mcX \sim \msP_1|_\mcX$.
In this case, however, $\msP_\frac{1}{4} \sim \msP_\frac{3}{4}$.
\end{remark}   


\subsection{Pairwise orthogonality implies  
         orthogonality in aggregate}
\label{subsec:PairwiseOrthogonality}

In this subsection, all measures $\msP_g$, $\msP_\mu$ are
considered as being defined on the measurable space 
$(\widetilde{\mcA},\widetilde{\mcX})$
(see subsection \ref{subsec:TailSpace}).
We shall, however, drop the $\sim$ sign in subsequent notation.

We consider a mapping $h : \mcX \rightarrow \mcB(Q)$
defined as:
\begin{equation}   
\label{eq:hDef}
hA = \{ g \in Q \mid \msP_g(A)=1 \}
\end{equation}   

According to Kolmogorov's zero-one law, $\msP_g$ takes only values
$0$ and $1$; thus, $g \not\in hA$ is equivalent to $\msP_g(A)=0$.
The fact that $hA$ is a Borel set follows from proposition
\ref{pr:Meas2}, as $hA$ is the inverse image of Borel set $\{1\}$
under mapping $g \mapsto \msP_g(A)$.

Our first goal in this subsection is to show that $h$ is$\pmod{0}$
epimorphism of $\sigma$-algebras.

\begin{lemma}   
\label{lm:hIsHomo}
Mapping $h$ is a homomorphism of $\sigma$-algebras.
\end{lemma}

\begin{proof}
All required relations $h(A' \cap A'') = hA' \cap hA''$, etc.
and $h(\bigcap_n A_n) = \bigcap_n hA_n$, etc. can be
verified directly.
For example:
\begin{multline*}
{\textstyle
h(\bigcap_n A_n) = \{ g \in Q \mid \msP_g(\bigcap_n A_n)=1 \} =
\{ g \in Q \mid \forall n \msP_g(A_n)=1 \} = }
\\
{\textstyle
\bigcap_n \{ g \in Q \mid \msP_g(A_n)=1 \} = \bigcap_n hA_n }
\mbox{\qedhere}
\end{multline*}
\end{proof}   

A proof that $h$ is ``onto''$\pmod{0}$ requires establishing
additional facts.

\medskip

We shall extensively use the following theorem,
proved in \citet[][Theorem 2; the wording below is adopted for
our case]{Robbins:1948}:

\newtheorem*{RT}{Robbins' Theorem}   
\begin{RT}
Let $f : \mcA \rightarrow \mbR$ be an $\mcX$-measurable
nonnegative function.
Then the function $\hat{f}(g) \bydef \int_\mcA f(a) \,\msP_g(da)$
is $\mcB(Q)$-measurable and for every measure $\mu$ on $Q$
\begin{equation}   
\label{eq:Robbins}
\int_\mcA f(a) \,\msP_\mu(da) = \int_Q \hat{f}(g) \,\mu(dg) =
\int_Q \left( \int_\mcA f(a) \,\msP_g(da) \right) \,\mu(dg)
\end{equation}   
\end{RT}   

The first important corollary of the Robbins' theorem is:

\begin{lemma}   
\label{lm:hPreservesMeasure}
For every measure $\mu$ on $Q$ and for every $A \in \mcX$,
$\msP_\mu(A) = \mu(hA)$.
\end{lemma}

\begin{proof}%
\begin{multline*}   
\msP_\mu(A) = \int_\mcA I_A(a) \,\msP_\mu(da) =
\int_Q \left( \int_\mcA I_A(a) \,\msP_g(da) \right) \,\mu(dg) =\\
\int_Q \msP_g(A) \,\mu(dg) = \int_Q I_{hA}(g) \,\mu(dg) = \mu(hA)
\mbox{\qedhere}
\end{multline*}   
\end{proof}   

\begin{lemma}   
\label{lm:AbsContinuity}
$\forall \mu', \mu'' \in \mcP(Q) :
\mu' \ll \mu'' \impl \msP_{\mu'} \ll \msP_{\mu''}$.
\end{lemma}

\begin{proof}%
$\msP_{\mu''}(A)=0 \eqval \mu''(hA)=0 \impl 
\mu'(hA)=0 \eqval \msP_{\mu'}(A)=0$. 
\end{proof}   

\begin{lemma}   
\label{lm:EquivalentMeasure}
Let $\mu', \mu'' \in \mcP(Q)$ be such that
$\msP_{\mu'} \ll \msP_{\mu''}$.
Then there exists a measure $\bar{\mu}'' \in \mcP(Q)$ such that
$\bar{\mu}'' \ll \mu''$ and $\msP_{\bar{\mu}''} = \msP_{\mu'}$.
\end{lemma}

\begin{proof}
Let
\begin{equation*}   
z(a) = \frac{d\msP_{\mu'}}{d\msP_{\mu''}}(a)
\end{equation*}   
be a Radon-Nikodym derivative of $\msP_{\mu'}$ with respect to
$\msP_{\mu''}$.
Let us fix any variant of $z(a)$.
By Robbins' theorem, the function 
$y(g) = \int_\mcA z(a) \,\msP_g(da)$ is $\mcB(Q)$-measurable
and
\begin{equation*}   
\int_Q y(g) \,\mu''(dg) =
\int_Q \left( \int_\mcA z(a) \,\msP_g(da) \right) \,\mu''(dg) =
\int_\mcA z(a) \,\msP_{\mu''}(da) = 1
\end{equation*}   
Thus, the set function $\bar{\mu}''(B) \bydef \int_B y(g) \,\mu''(dg)$
is a probabilistic measure on $Q$ and $\bar{\mu}'' \ll \mu''$.
Further, using lemma \ref{lm:hPreservesMeasure} and Robbins' theorem,
one obtains for arbitrary $A \in \mcX$:
\begin{multline*}   
\msP_{\bar{\mu}''}(A) = \bar{\mu}''(hA) =
\int_Q I_{hA}(g) y(g) \,\mu''(dg) = \\
\int_Q I_{hA}(g) \left( \int_\mcA z(a) \,\msP_g(da) \right)
                                          \,\mu''(dg) \overset{(*)}{=}
\int_Q \left( \int_\mcA I_A(a) z(a) \,\msP_g(da) \right) \,\mu''(dg) =\\
\int_\mcA I_A(a) z(a) \,\msP_{\mu''}(da) = \msP_{\mu'}(A)
\end{multline*}   
The intermediate equality $(*)$ is verified directly:
\begin{equation*}   
\msP_g(A)=1 \impl I_{hA}(g)=1 \text{ and } 
\int_A z(a) \,\msP_g(da) = \int_\mcA z(a) \,\msP_g(da)
\end{equation*}   
and
\begin{equation*}   
\msP_g(A)=0 \impl I_{hA}(g)=0 \text{ and }
\int_A z(a) \,\msP_g(da) = 0 
\end{equation*}   

Thus, one obtains $\msP_{\bar{\mu}''} = \msP_{\mu'}$, q.e.d.
We did not prove that different variants of $z(a)$ lead to the
same measure $\bar{\mu}''$, but it is insignificant for this lemma,
as it claims only {\em existence} of measure $\bar{\mu}''$ with
desired properties, and any variant of $z(a)$ gives such example.
\end{proof}   

\begin{lemma}   
\label{lm:noAtoms}
Suppose that 
$\forall g', g'' \in Q : g' \neq g'' \impl \msP_{g'} \perp \msP_{g''}$.
Then continuity of $\mu \in \mcP(Q)$ implies continuity of $\msP_\mu$.
\end{lemma}

\begin{proof}
Assume, in contrary, that $A_0 \in \mcX$ is an atom of $\msP_\mu$,
i.e. $\msP_\mu(A_0) > 0$ and for every $A \in \mcX$, $A \subseteq A_0$,
either $\msP_\mu(A) = \msP_\mu(A_0)$ or $\msP_\mu(A) = 0$.

Let $B_0 = hA_0$. Then $\mu(B_0)>0$, and we can consider
a measure $\mu_0 = \mu(\cdot|B_0)$.
It is easy to see that $\msP_{\mu_0}$ is carried by $A_0$ and
$A_0$ is an atom of $\msP_{\mu_0}$.

Consider any measure $\mu' \ll \mu_0$. By lemma \ref{lm:AbsContinuity},
$\msP_{\mu'} \ll \msP_{\mu_0}$. Thus, there exists a Radon-Nikodym
derivative $z(a)=d\msP_{\mu'}/d\msP_{\mu_0}$.
But as $A_0$ is atom of $\msP_{\mu_0}$, $z(a)$ is $\msP_{\mu_0}$-a.s.
a constant on $A_0$, and it follows from normalization conditions
that this constant is $1$. 
Thus, $\msP_{\mu'}=\msP_{\mu_0}$ for every measure $\mu' \ll \mu_0$.

Now return to the set $B_0$. As $\mu_0(B_0)=1$ and $\mu_0$ has
no atoms, it is possible to construct two sequences of subsets
of $B_0$, $\{B'_n\}_n$ and $\{B''_n\}_n$, such that
$B'_n \downdownarrows \{g'\}$, $B''_n \downdownarrows \{g''\}$,
$g' \neq g''$, and $\mu_0(B'_n)>0$, $\mu_0(B''_n)>0$ for all $n$.

Consider measures $\mu'_n = \mu_0(\cdot|B'_n)$ and
$\mu''_n = \mu_0(\cdot|B''_n)$.
One has $\mu'_n \xrightarrow{w} \delta_{g'}$ and
$\mu''_n \xrightarrow{w} \delta_{g''}$,
while $\msP_{\mu'_n} = \msP_{\mu''_n} = \msP_{\mu_0}$ for all $n$.
From continuity of the mixing operator \citep{Kovtun:2005d},
one obtains
\begin{equation*}   
\msP_{g'} = \lim_{n \rightarrow \infty} \msP_{\mu'_n} =
\lim_{n \rightarrow \infty} \msP_{\mu''_n} = \msP_{g''}
\end{equation*}   
which contradicts to the assumption of the lemma that 
$\msP_{g'} \perp \msP_{g''}$.
\end{proof}   

\begin{lemma}   
\label{lm:ContOrtDelta}
Suppose that 
$\forall g', g'' \in Q : g' \neq g'' \impl \msP_{g'} \perp \msP_{g''}$.
Let $\mu \in \mcP(Q)$ be a continuous measure.
Then $\msP_\mu \perp \msP_g$ for every $g \in Q$.
\end{lemma}

\begin{proof}
Assume premise of the lemma and take arbitrary continuous $\mu$
and arbitrary $g$.
Let $\bar{\mu} = \fot(\mu + \delta_g)$.
Due to linearity of the mixing operator,
$\msP_{\bar{\mu}} = \fot(\msP_\mu + \msP_g)$,
and thus $\msP_\mu\ll\msP_{\bar{\mu}}$ and $\msP_g\ll\msP_{\bar{\mu}}$.
Let
\begin{equation*}   
\bar{z}(a)=\frac{ d\msP_\mu }{ d\msP_{\bar{\mu}} }(a), \qquad
z(a) = \frac{ d\msP_g }{ d\msP_{\bar{\mu}} }(a), \qquad
A_0 = \{ a \mid \bar{z}(a) \cdot z(a) > 0 \}
\end{equation*}   
Then $\msP_\mu \perp \msP_g$ is equivalent to $\msP_{\bar{\mu}}(A_0)=0$.

Assume in contrary that $\msP_{\bar{\mu}}(A_0)>0$.
Take arbitrary $A_1 \subseteq A_0$.
If $0 < \msP_{\bar{\mu}}(A_1) < \msP_{\bar{\mu}}(A_0)$,
then $\msP_g(A_1)>0$ and $\msP_g(A_0 \setminus A_1)>0$,
which contradicts to the Kolmogorov's zero-or-one law.
Thus, $A_0$ is an atom of $\msP_{\bar{\mu}}$, and consequently
an atom of $\msP_\mu$, which contradicts to the lemma
\ref{lm:noAtoms}.
\end{proof}   

\begin{lemma}   
\label{lm:ContOrtCount}
Suppose that 
$\forall g', g'' \in Q : g' \neq g'' \impl \msP_{g'} \perp \msP_{g''}$.
Let $\mu \in \mcP(Q)$ be a continuous measure,
and $\mu' \in \mcP(Q)$ be a counting measure.
Then $\msP_\mu \perp \msP_{\mu'}$.
\end{lemma}

\begin{proof}
Any counting measure can be represented in form
$\mu' = \sum_{i=1}^\infty a_i \delta_{g_i}$,
where $a_i \ge 0$ and $\sum a_i = 1$.
Then due to linearity of the mixing operator
$\msP_{\mu'} = \sum a_i \msP_{g_i}$.

For every $i=1,\dots$ take $A_i$ such that $\msP_\mu(A_i)=1$
and $\msP_{g_i}(A_i)=0$ (such sets exist due to lemma
\ref{lm:ContOrtDelta}).
Now take $A = \bigcap_i A_i$.
Then $\msP_\mu(A)=1$ and $\msP_{\mu'}(A)=0$.
\end{proof}   

\begin{lemma}   
\label{lm:hOntoClosed}
Suppose that 
$\forall g', g'' \in Q : g' \neq g'' \impl \msP_{g'} \perp \msP_{g''}$.
Let $\mu \in \mcP(Q)$ be a continuous measure, $F$ be a closed
set such that $\mu(F)>0$, $\mu' = \mu(\cdot|F)$,
and $z(a)$ be a Radon-Nikodym derivative
of $\msP_{\mu'}$ with respect to $\msP_\mu$.
Let $A_0 = \{ a \mid z(a)>0 \}$.
Then $\mu(F \bigtriangleup hA_0)=0$.
\end{lemma}

\begin{proof}
As $\mu'(hA_0) = \msP_{\mu'}(A_0) = 1$, 
one obtains $\mu'(F \setminus hA_0) = 0$; consequently,
$\mu(F \setminus hA_0) = 0$.

Next, assume in contrary that $\mu(hA_0 \setminus F) > 0$.
As $F$ is closed, one can construct a sequence of
sets $\{B_n\}_n$ such that $B''_n \subseteq hA_0 \setminus F$, 
$\mu(B_n)>0$,
and $B_n \downdownarrows \{g\}$ for some $g \in hA_0 \setminus F$.
Let $\mu_n = \mu(\cdot|B_n)$. Then $\mu_n \ll \mu$,
and by lemma \ref{lm:AbsContinuity} $\msP_{\mu_n}\ll\msP_\mu$.
Let
\begin{equation*}   
z_n(a) = \frac{d\msP_{\mu_n}}{d\msP_\mu}(a), \qquad
A_n = \{ a \mid z_n(a) > 0 \}
\end{equation*}   
It is always possible to choose such variants of $z_n(a)$
that $A_n \subseteq A_0$, $A_{n+1} \subseteq A_n$.
Let $B'_n = hA_n \cap F$. Then $B'_{n+1} \subseteq B'_n$
and $\mu'(B'_n)=\mu'(hA_n)=\msP_{\mu'}(A_n)>0$
(as $\msP_{\mu'}(A_n)=0 \eqval \msP_\mu(A_n)=0$).

Let $\mu'_n = \mu'(\cdot|B'_n)$.
We claim that $\msP_{\mu'_n} \ll \msP_{\mu_n}$.
To prove this, take any set $A$ such that $\msP_{\mu'_n}(A)>0$;
then $\msP_{\mu'_n}(A) = \mu'_n(hA) = \mu'_n(hA \cap hA_n) =
\msP_{\mu'_n}(A \cap A_n)$; further, as $\mu'_n \ll \mu' \ll \mu$,
$\msP_\mu(A \cap A_n)>0$, and consequently 
$\msP_{\mu_n}(A \cap A_n)>0$, q.e.d.

Further, we claim that there exists a point $g' \in F$
such that $\forall\,\varepsilon>0 \,\,\forall\,n : 
\mu'_n(U_\varepsilon g')>0$.
To prove this, suppose the contrary, i.e. that for every
point $g \in F$ there exists $\varepsilon(g)$ such that
$\mu'_n(U_{\varepsilon(g)}g)=0$ for all $n$ starting some $n(g)$
(note that for $n'>n$ one has $\mu'_{n'}\ll\mu'_n$, and thus
$\mu'_n(U)=0 \impl \mu'_{n'}(U)=0$).
This system covers $F$, and as $F$ is compact, it contains
a finite subcover $\{U_{\varepsilon(g_i)}g_i\}_{i=1}^m$.
Let $n_0 = \max(n(g_1),\dots,n(g_m))$.
Then 
$1 = \mu'_{n_0}(F) \le \sum_i \mu'_{n_0}(U_{\varepsilon(g_i)}g_i) = 0$.
This contradiction proves the claim.

Now take $B''_n = B'_n \cap U_{1/n}g' \cup \{g'\}$
and $\mu''_n = \mu'_n(\cdot|B''_n)$.
One has $B''_n \downdownarrows \{g'\}$. 
As $\mu''_n \ll \mu'_n$, one obtains 
$\msP_{\mu''_n} \ll \msP_{\mu'_n} \ll \msP_{\mu_n}$.
By lemma \ref{lm:EquivalentMeasure}, one can construct measures
$\bar{\mu}_n \ll \mu_n$ such that $\msP_{\bar{\mu}_n} = \msP_{\mu''_n}$.
Finally, using continuity of the mixing operator \citep{Kovtun:2005d},
one obtains
\begin{equation*}   
\msP_{g'} = \lim_{n \rightarrow \infty} \msP_{\mu''_n} =
\lim_{n \rightarrow \infty} \msP_{\bar{\mu}_n} = \msP_{g}
\end{equation*}   
which contradicts to the assumption of the lemma that 
$\msP_{g'} \perp \msP_{g}$.
\end{proof}   

\begin{lemma}   
\label{lm:hOntoCont}
Suppose that 
$\forall g', g'' \in Q : g' \neq g'' \impl \msP_{g'} \perp \msP_{g''}$.
Then for every continuous measure $\mu \in \mcP(Q)$
and for every $B \in \mcB(Q)$ 
there exists $A \in \mcX$ such that $\mu(B \bigtriangleup hA)=0$.
\end{lemma}

\begin{proof}
Take arbitrary continuous $\mu \in \mcP(Q)$.
Let $\mcB_0 = \{ B \in \mcB(Q) \mid 
\exists\,A \in \mcX : \mu(B \bigtriangleup hA)=0 \}$.
By lemma \ref{lm:hOntoClosed}, all closed subsets of $Q$
belong to $\mcB_0$.
As $h$ is a homomorphism of $\sigma$-algebras (lemma \ref{lm:hIsHomo}),
$\mcB_0$ is closed under finite and countable unions and intersections
and under complements. Thus, $\mcB_0$ coincides with $\mcB(Q)$.
\end{proof}   

\begin{lemma}   
\label{lm:hOntoCount}
Suppose that 
$\forall g', g'' \in Q : g' \neq g'' \impl \msP_{g'} \perp \msP_{g''}$.
Then for every counting measure $\mu \in \mcP(Q)$
and for every $B \in \mcB(Q)$ 
there exists $A \in \mcX$ such that $\mu(B \bigtriangleup hA)=0$.
\end{lemma}

\begin{proof}
Let $\mu = \sum_{i=1}^\infty a_i \delta_{g_i}$,
where $a_i \ge 0$ and $\sum a_i = 1$.
Then $\msP_\mu = \sum a_i \msP_{g_i}$.

For every $n \neq m$ take $A_{nm}$ such that $\msP_{g_n}(A_{nm})=1$
and $\msP_{g_m}(A_{nm})=0$, and take $A_n = \bigcap_{m \neq n} A_{nm}$.
Then $\msP_{g_n}(A_n)=1$ and for every $m \neq n$ $\msP_{g_m}(A_n)=0$.

Now for every $B \in \mcB(Q)$ take 
$A_B = \bigcup \{ A_n \mid g_n \in B \}$.
It is easy to verify that $\mu(B \bigtriangleup hA_B) = 0$.
\end{proof}   

\begin{lemma}   
\label{lm:hOnto}
Suppose that 
$\forall g', g'' \in Q : g' \neq g'' \impl \msP_{g'} \perp \msP_{g''}$.
Then for every measure $\mu \in \mcP(Q)$
and for every $B \in \mcB(Q)$ 
there exists $A \in \mcX$ such that $\mu(B \bigtriangleup hA)=0$.
\end{lemma}

\begin{proof}
Every measure $\mu$ on $Q$ can be represented as
$\alpha \mu' + (1-\alpha)\mu''$, where $\mu'$ is continuous
and $\mu''$ is counting.

By lemma \ref{lm:ContOrtCount}, one can find sets $A',A'' \in \mcX$
such that $\msP_{\mu'}(A')=1$, $\msP_{\mu''}(A'')=1$,
and $A' \cap A'' = \varnothing$.

Now take arbitrary $B \in \mcB(Q)$. By lemma \ref{lm:hOntoCont},
there exists $A_1 \in \mcX$ such that $\mu'(B \bigtriangleup hA_1)=0$,
and by lemma \ref{lm:hOntoCount}, there exists $A_2 \in \mcX$
such that $\mu''(B \bigtriangleup hA_2)=0$.
It is always possible to choose these sets such that
$A_1 \subseteq A'$ and $A_2 \subseteq A''$.
Thus, $A_1 \cap A_2 = \varnothing$ and $\mu'(hA_2)=\mu''(hA_1)=0$.

Finally,
\begin{multline*}
\mu(B \bigtriangleup h(A_1 \cup A_2)) =
\mu(B \bigtriangleup hA_1 \bigtriangleup hA_2) = \\
\alpha \mu'(B \bigtriangleup hA_1 \bigtriangleup hA_2) +
(1-\alpha) \mu''(B \bigtriangleup hA_1 \bigtriangleup hA_2) = \\
\alpha \mu'(B \bigtriangleup hA_1) +
(1-\alpha) \mu''(B \bigtriangleup hA_2) =
\alpha \cdot 0 + (1-\alpha) \cdot 0 = 0
\mbox{\qedhere}
\end{multline*}
\end{proof}   

Lemma \ref{lm:hOnto} allows us to obtain a corollary
which is of interest by itself.
Namely:

\begin{lemma}   
\label{lm:OrtMixIsOrt}
Suppose that 
$\forall g', g'' \in Q : g' \neq g'' \impl \msP_{g'} \perp \msP_{g''}$.
Then $\forall \, \mu_1, \mu_2 \in \mcP(Q) :
\mu_1 \perp \mu_2 \impl \msP_{\mu_1} \perp \msP_{\mu_2}$.
\end{lemma}

\begin{proof}
Assume premise of the lemma and take arbitrary $\mu_1$, $\mu_2$
such that $\mu_1 \perp \mu_2$.
Then there exists Borel $B_1, B_2 \subseteq Q$ such that
$\mu_1(B_1)=1$, $\mu_2(B_2)=1$, and $B_1 \cap B_2 = \varnothing$.

Let $\mu = \fot(\mu_1 + \mu_2)$.
By lemma \ref{lm:hOnto} there exist sets $A'_1, A'_2 \in \mcX$ such that
$\mu(B_1 \bigtriangleup hA'_1) = \mu(B_2 \bigtriangleup hA'_2) = 0$.
As $\msP_\mu(A'_1 \cap A'_2)=\mu(hA'_1 \cap hA'_2)=\mu(B_1 \cap B_2)=0$,
the sets $A_1 = A'_1 \setminus (A'_1 \cap A'_2)$ and
$A_2 = A'_2 \setminus (A'_1 \cap A'_2)$ also satisfy
$\mu(B_1 \bigtriangleup hA_1) = \mu(B_2 \bigtriangleup hA_2) = 0$.

Now one has $\msP_{\mu_1}(A_1) = \mu_1(hA_1) = \mu_1(B_1) = 1$;
similarly, $\msP_{\mu_2}(A_2) = 1$.
As $A_1$ and $A_2$ are disjoint, one obtains orthogonality of
$\msP_{\mu_1}$ and $\msP_{\mu_2}$.
\end{proof}   

On the other hand, lemma \ref{lm:hOnto} can be easily derived
from lemma \ref{lm:OrtMixIsOrt}.
For this, take arbitrary $\mu \in \mcP(Q)$ and $B \in \mcB(Q)$.
If $\mu(B)=1$ or $\mu(B)=0$, one can take $A=\mcA$ or $A=\varnothing$,
respectively.
Otherwise, one can consider measures $\mu(\cdot|B)$ and
$\mu(\cdot|Q \setminus B)$.
These measures are orthogonal, and by lemma \ref{lm:OrtMixIsOrt}
measures $\msP_{\mu(\cdot|B)}$ and $\msP_{\mu(\cdot|Q \setminus B)}$
are orthogonal as well.
Take $A$ such that $\msP_{\mu(\cdot|B)}(A)=1$ and
$\msP_{\mu(\cdot|Q \setminus B)}(A)=0$.
It is easy to see that $\mu(B \bigtriangleup hA)=0$, which proves
lemma \ref{lm:hOnto}.

\medskip

Lemma \ref{lm:hOnto} shows that for every $\mu \in \mcP(Q)$
the metric structure of spaces $(Q,\mcB(Q),\mu)$ and
$(\mcA,\mcX,\msP_\mu)$ are isomorphic.
This allows us to employ Rokhlin's technique
\citep[see][]{Rokhlin:1949,Rokhlin:1952}
to obtain the main result of the present subsection.

As $Q$ is a complete separable metric space,
$(Q,\mcB_\mu(Q),\bar{\mu})$ is a Lebesgue space for every
Borel measure $\mu$
(recall that $\mcB_\mu(Q)$ denotes a Lebesgue completion of Borel
$\sigma$-algebra $\mcB(Q)$ with respect to measure $\mu$,
and $\bar{\mu}$ denotes the corresponding completion of measure $\mu$).

Let us fix some measure $\mu$ on $(Q,\mcB(Q))$, and
let $\mfG=\{\Gamma_i\}_{i=1}^\infty$ be a topology base of $Q$;
it is also a Rokhlin's basis of $(Q,\mcB_\mu(Q),\bar{\mu})$.
Let further $\hat{\Gamma}_i$ be a coset of $\Gamma_i$
by ideal of $\bar{\mu}$-negligible sets $I_{\bar{\mu}}$.

By lemma \ref{lm:hOnto}, there exists a set $C_i \in \mcX$
such that $hC_i \in \hat{\Gamma}_i$.
Let $\Delta_i = \Gamma_i \bigtriangleup hC_i$,
and let $Q_0 = \bigcup_{i=1}^\infty \Delta_i$,
$Q_1 = Q \setminus Q_0$.
Then we have exact equalities $hC_i = \Gamma_i \cap Q_1$
and $h(\mcA \setminus C_i) = Q_1 \setminus \Gamma_i$.

Note that both $Q_0$ and $Q_1$ are Borel sets
(as all $\Gamma_i$ are open, and all $hC_i$ are Borel).

For arbitrary $g \in Q_1$, let $\Upsilon_i(g) = \Gamma_i \cap Q_1$
if $g \in \Gamma_i$, and $\Upsilon_i(g) = Q_1 \setminus \Gamma_i$
if $g \not\in \Gamma_i$.
Correspondingly, let $D_i(G)=C_i$ if $g \in \Gamma_i$, and
$D_i(g)=\mcA \setminus C_i$ if $g \not\in \Gamma_i$.
Then $hD_i(g) = \Upsilon_i(g)$.

Finally, let $A_g = \bigcap_i D_i(g)$.
As $h$ is homomorphism of $\sigma$-algebras, one obtains
$hA_g = h(\bigcap_i D_i(g)) = \bigcap_i hD_i(g) =
\bigcap_i \Upsilon_i(g) = \{g\}$.
This means that $\msP_g(A_g)=1$ and for every $g' \neq g$,
$\msP_{g'}(A_g)=0$.
Thus, we proved:

\begin{lemma}   
\label{lm:SingSets}
Suppose that 
$\forall g', g'' \in Q : g' \neq g'' \impl \msP_{g'} \perp \msP_{g''}$.
Then for every measure $\mu \in \mcP(Q)$
there exists a Borel subset $Q_1 \subseteq Q$ of full measure,
$\mu(Q_1)=1$, 
and a mapping $\psi : Q_1 \rightarrow \mcX : g \mapsto A_g$ satisfying
(a) $\forall\, g \in Q_1 : \msP_g(A_g)=1$ and
(b) $\forall\,g,g'\in Q_1: g \neq g' \impl A_g \cap A_{g'}=\varnothing$.
\end{lemma}   

As sets $\{A_g\}_{g \in Q_1}$ are disjoint, one can consider
a mapping $\mcA \rightarrow Q : a \mapsto g_a$, defined as
$g_a = g$, if $a \in A_g$, and $g_a = g_0$ for 
$a \not\in \bigcup_{g \in Q_1} A_g$; $g_0 \in Q$ is taken
arbitrarily.
It is easy to see that the mapping $a \mapsto g_a$ is
$\mcB(Q)$-$\mcX$-measurable.
Thus, we have

\begin{lemma}   
\label{lm:InvSingSets}
Suppose that 
$\forall g', g'' \in Q : g' \neq g'' \impl \msP_{g'} \perp \msP_{g''}$.
Then for every measure $\mu \in \mcP(Q)$
there exists a Borel subset $Q_1 \subseteq Q$ of full measure,
$\mu(Q_1)=1$,
and a $\mcB(Q)$-$\mcX$-measurable mapping
$\varphi : \mcA \rightarrow Q : a \mapsto g_a$
satisfying $Q_1 \subseteq \varphi(\mcA)$ and
$\msP_g(\varphi^{-1}(g)) = 1$ for every $g \in Q_1$.
\end{lemma}   


\subsection{Conditional probabilities with respect  
         to tail $\sigma$-algebra}
\label{subsec:ConditionalProbabilities}

Let $\msP$ be a probabilistic measure on $(\mcA,\mcB(\mcA))$.
We say that $\msP$ {\em satisfies 0-1-law}, if for every
$A \in \mcX$ either $\msP(A)=1$ or $\msP(A)=0$.

Note that all measures $\msP_g$, $g \in Q$, satisfy 0-1-law.

For measures satisfying 0-1-law, conditional probabilities
with respect to tail $\sigma$-algebra $\mcX$ can be expressed
directly. Namely, a well-known fact is:

\begin{lemma}   
\label{lm:CondPr1}
Let measure $\msP$ satisfy 0-1-law.
Then for every $B \in \mcB(\mcA)$ one has
$\msP(B | \mcX)(a) = \msP(B)$ for $\msP$-almost all $a \in \mcA$. 
\end{lemma}


With this we can easily prove:

\begin{lemma}   
\label{lm:TailSufficient}
Suppose that 
$\forall g', g'' \in Q : g' \neq g'' \impl \msP_{g'} \perp \msP_{g''}$.
Then for every measure $\mu \in \mcP(Q)$ there exists a Borel
set $Q_1$ of full measure, $\mu(Q_1)=1$, such that
tail $\sigma$-algebra $\mcX$ is sufficient for the family
$\{\msP_g\}_{g \in Q_1}$.
\end{lemma}

\begin{proof}
Take $Q_1$ satisfying conditions of lemma \ref{lm:InvSingSets}.
For $a \in \mcA$ and $B \in \mcB(\mcA)$, let
$P(a,B) \bydef \msP_{\varphi(a)}(B)$
(where $\varphi$ was defined in lemma \ref{lm:InvSingSets}).
By definition, $P(a,B)$ is a probabilistic measure on
$(\mcA,\mcB(\mcA))$ for every fixed $a \in \mcA$.
Further, by proposition \ref{pr:Meas2} and lemma \ref{lm:InvSingSets},
$P(a,B)$ is $\mcX$-measurable as function of $a$ for every
fixed $B$,
and by lemma \ref{lm:InvSingSets}, for every $g \in Q_1$
and for every fixed $B$,
$P(a,B) = \msP_g(B)$ for $\msP_g$-almost all $a$.
Finally, by lemma \ref{lm:CondPr1}, $P(a,B) = \msP_g(B | \mcX)(a)$
$\msP_g$-a.s., which completes the proof.
\end{proof}   

The inverse statement is also true:

\begin{lemma}   
\label{lm:InvTailSufficient}
Suppose that for every measure $\mu \in \mcP(Q)$ there exists a Borel
set $Q_1$ of full measure, $\mu(Q_1)=1$, such that
tail $\sigma$-algebra $\mcX$ is sufficient for the family
$\{\msP_g\}_{g \in Q_1}$.
Then
$\forall g', g'' \in Q : g' \neq g'' \impl \msP_{g'} \perp \msP_{g''}$.
\end{lemma}

\begin{proof}
Take arbitrary $g', g'' \in Q$ such that $g' \neq g''$
and consider a measure $\mu = \fot (\delta_{g'} + \delta_{g''})$.
Then any set $Q_1$ of full measure $\mu$ contains $g'$ and $g''$.
Sufficiency of $\mcX$ means that there exists a function
$P(a,B)$ such that $\msP_{g'}(B|\mcX)(a)=P(a,B)$ $\msP_{g'}$-a.s.
and $\msP_{g''}(B|\mcX)(a)=P(a,B)$ $\msP_{g''}$-a.s.
Family $\{\msP_g\}_{g \in Q}$ consists of distinct measures;
thus, there exists a set $B_0 \in \mcB(\mcA)$ such that
$\msP_{g'}(B_0) \neq \msP_{g''}(B_0)$.

Take $A' = \{ a \mid P(a,B_0)=\msP_{g'}(B_0) \}$ and
$A'' = \{ a \mid P(a,B_0)=\msP_{g''}(B_0) \}$.
By construction, $A' \cap A'' = \varnothing$, and
by lemma \ref{lm:CondPr1}, $\msP_{g'}(A')=\msP_{g''}(A'')=1$,
which proves orthogonality of $\msP_{g'}$ and $\msP_{g''}$.
\end{proof}   



\section{Discussion} 
\label{sec:Discussion}

\subsection{Orthogonality in aggregate: 
            for ``all'' or ``almost all'' $g \in Q$?}
\label{subsec:OrthogonalityInAggregate}

The result proven in lemma \ref{lm:SingSets} can be reformulated as
``pairwise orthogonality of measures $\{\msP_g\}_{g \in Q}$
implies orthogonality in aggregate for almost all measures
from this family.'' The ``almost all'' clause of the above statement
can be based on any measure $\mu$ on $Q$---note, however, that
sets $A_g$ may depend on the choice of $\mu$.
One can also see that there is no obvious modification of our proof,
which gives orthogonality in aggregate for the whole family
$\{\msP_g\}_{g \in Q}$.

A natural question arise: whether the absence of full orthogonality
in aggregate is an immanent property or it is merely a deficiency
of our proof?

This question can be also formulated in terms of lemma \ref{lm:hOnto},
which states that homomorphism $h : \mcX \rightarrow \mcB(Q)$ is
``onto''$\pmod{0}$; again, ``$\!\!\pmod{0}$'' holds for any measure
$\mu$ on $Q$---and again, set $A$ satisfying $\mu(B \bigtriangleup A)=0$
may depend on the choice of $\mu$.
This point of view can be (less formally) expressed as:
is the tail $\sigma$-algebra $\mcX$ sufficiently rich
to have a factor-algebra isomorphic to $\mcB(Q)$,
or one can obtain only isomorphism$\pmod{0}$?

\subsection{Structure of tail $\sigma$-algebras} 
\label{subsec:StructureOfTail}

Structure of tail $\sigma$-algebras deserved not too much attention.
Such lack of interest can be explained,
first, by the fact that tail $\sigma$-algebras are ``trivial''
from probabilistic point of view (any independent distribution,
being restricted to the tail $\sigma$-algebra,
takes only values 0 and 1),
and second, by the absence of topology or other ``good'' structure
generating tail $\sigma$-algebras (cf., for example,
a comprehensive study of Borel sets by means of topology and
descriptive set theory).

A study performed in the present paper demonstrates that an explicit
investigation of properties of tail $\sigma$-algebra might
bring fruitful results:
technical basis for our main theorem is existence of homomorphism
$h : \mcX \rightarrow \mcB(Q)$, which is epimorphism$\pmod{0}$,
proved in lemma \ref{lm:hOnto}.
An important corollary of the latest fact is that for every
measure $\mu$ on $(Q,\mcB(Q))$ the metric structures of spaces
$(Q,\mcB(Q),\mu)$ and $(\mcA,\mcX,\msP_\mu)$ are isomorphic.

Note that $\mcX$ is constructed from ranges of random variables
$X_1,\dots$; nothing from distribution of $X_1,\dots$ is involved
in construction of $\mcX$.
Lemma \ref{lm:hOnto} establishes existence of homomorphism $h$
for {\em any} set $Q$---provided that $Q$ is compact and continuously
parameterize pairwise orthogonal family of independent measures.
In other words, $\mcX$ is universal with respect to family
of $\sigma$-algebras $\{\mcB(Q)\}_Q$.
It looks like this universality is a characteristic property of tail
$\sigma$-algebras, and might be a useful tool for further
investigation of their properties.

\subsection{Connection with the strong law  
            of large numbers}
\label{subsec:ConnectionWithTheStrong}

One fact that we have proved in this article is that
pairwise orthogonality of considered family of independent
distributions implies their orthogonality in aggregate.
Informally speaking, one can say the sets $A_g$, provided by
lemma \ref{lm:SingSets} are ``sufficiently small'' to be
disjoint.

Another way to obtain ``small'' sets with $\msP_g$ measure 1
is provided by the strong law of large numbers.
In some sense, the ``smallness'' of the set of outcomes
provided by the strong law of large numbers is one of
cornerstones of probability theory.
When random variables $X_j$ are identically distributed
for every $g$, the sets given by the law of large numbers
for different $\msP_g$ are disjoint.

Consider, however, the following example.

Let $X_j$ be a sequence of binary random variables
taking values 0 and 1, $Q = [-1,1]$, and for every $g \in Q$,
$\msP_g(X_j=1)=\frac{1}{2}+\frac{g}{2\sqrt{j}}$ and
$\msP_g(X_j=0)=\frac{1}{2}-\frac{g}{2\sqrt{j}}$.
It is easy to verify that condition (\ref{mti:Hellinger2})
of the main theorem is satisfied; thus, distributions $\msP_g$
are pairwise orthogonal.
Now from lemma \ref{lm:SingSets} one obtains existence
of a disjoint family of sets $\{A_g\}_{g \in Q_1}$
(where $Q_1$ has Lebesgue measure 1, and thus has cardinality
of continuum) such that $\msP_g(A_g)=1$.

The strong law of large numbers in form of Kolmogorov
\citep[][IV.3.2]{Shiryaev:2004}, adopted for our case,
states that for every increasing unbounded sequence of
positive reals $\{b_n\}_n$ such that
$\sum_1^\infty \frac{1}{b_n^2} < \infty$ the set
\begin{equation*}
A_g^{(LN)} =
\{ a \in \mcA \mid
    {\textstyle \frac{\sum_1^n a_j - \sum_1^n \msP_g(X_j=1)}{b_n} }
        \rightarrow 0 \}
\end{equation*}
has $\msP_g$ probability 1.

Using the facts that $\sum_1^\infty \frac{1}{b_n^2} < \infty$
implies that $\frac{\sqrt{n}}{b_n} \rightarrow 0$
and that
\begin{equation*}
\sum_1^n \msP_g(X_j=1) - \sum_1^n \msP_{g'}(X_j=1) =
\sum_1^n \frac{g-g'}{2\sqrt{j}} = O(\sqrt{n})
\end{equation*}
one obtains that $A_g^{(LN)} = A_{g'}^{(LN)}$
for all $g,g' \in Q$.

Thus, the strong law of large numbers may produce sets
of probability 1 that can be made
``significantly smaller''---namely, the set produced
by the strong law of large numbers can be splitted into
continuum of subsets, each of which has probability 1
with respect to corresponding distribution.

We believe that further investigation of consequences of
this fact might bring interesting results.

\subsection{Constructiveness of proofs} 
\label{subsec:ConstructivenessOfProofs}

The proof of existence of set $A_g$ in lemma \ref{lm:SingSets}
is non-constructive; moreover, it explicitly uses axiom of choice.
It would be desirable to find a constructive description
of sets $A_g$ (e.g., as in the strong law of large numbers),
as it would simplify analysis of their properties.
In particular, we expect that constructive description of sets
$A_g$ would clarify the answer to the question of section
\ref{subsec:OrthogonalityInAggregate}.


\section{Acknowledgements} 
\label{sec:Acknowledgements}

The first author personally thanks professor James N. Siedow
(Vice Provost for Research, Duke University),
whose trust in the author made this work possible.


\bibliographystyle{biostatistics}
\bibliography{math,Probability,Mixtures,LSA,GoM,mkovtun}

\end{document}